\documentstyle[11pt,fullpage]{article}

\newtheorem{thm}{Theorem}[section]
\newtheorem{prop}[thm]{Proposition}
\newtheorem{lem}[thm]{Lemma}
\newtheorem{cor}[thm]{Corollary}
\newtheorem{obs}[thm]{Observation}
\newtheorem{defn}[thm]{Definition}

\newtheorem{expl}[thm]{Example}
\newtheorem{rem}[thm]{Remark}

 \def\pushright#1{{
    \parfillskip=0pt            
    \widowpenalty=10000         
    \displaywidowpenalty=10000  
    \finalhyphendemerits=0      
    \leavevmode                 
    \unskip                     
    \nobreak                    
    \hfil                       
    \penalty50                  
    \hskip.2em                  
    \null                       
    \hfill                      
    {#1}                        
    \par}}                      
 \def\qed{\pushright{\rule{2mm}{3mm}}\penalty-700 \smallskip}

\newenvironment{proof}[1]{\begin{trivlist} \item[{\bf ~Proof}#1.]}%
{\qed\end{trivlist}}

\newenvironment{prf}[1]{\begin{trivlist} \item[{\bf ~Proof}#1.]}%
{\qed\end{trivlist}}

\newcounter{countroman}
{\begin{list}{{\rm (\roman{countroman})}}{\usecounter{countroman}}}%
{\end{list}}

\newcommand{\relcirc}{\mathrel{\circ}}   
\newcommand{\oto}{\mathrel{\relbar\joinrel\relcirc}}     
\newcommand{\x}{\times}

\newcommand{\ox}{\otimes}

\def\mathsl#1{\expandafter\def\csname#1\endcsname{\mathop{\sl #1}\nolimits}}
\def\mathbf#1{\expandafter\def\csname#1\endcsname{\mathop{\rm\bf #1}\nolimits}}

\def\mathopn#1{\expandafter\def\csname#1\endcsname{\mathop{\rm #1}\nolimits}}
\def\mathname#1{\expandafter\def\csname#1\endcsname{\mathop{\bf #1}\nolimits}}
\mathopn{op}
\mathopn{comp}

\makeatletter

\newdimen\w@dth

\def\setw@dth#1#2{\setbox\z@\hbox{\scriptsize $#1$}\w@dth=\wd\z@
\setbox\@ne\hbox{\scriptsize $#2$}\ifnum\w@dth<\wd\@ne \w@dth=\wd\@ne \fi
\advance\w@dth by 1.2em}

\def\t@^#1_#2{\allowbreak\def\n@one{#1}\def\n@two{#2}\mathrel
{\setw@dth{#1}{#2}
\mathop{\hbox to \w@dth{\rightarrowfill}}\limits
\ifx\n@one\empty\else ^{\box\z@}\fi
\ifx\n@two\empty\else _{\box\@ne}\fi}}
\def\t@@^#1{\@ifnextchar_ {\t@^{#1}}{\t@^{#1}_{}}}

\def\t@left^#1_#2{\def\n@one{#1}\def\n@two{#2}\mathrel{\setw@dth{#1}{#2}
\mathop{\hbox to \w@dth{\leftarrowfill}}\limits
\ifx\n@one\empty\else ^{\box\z@}\fi
\ifx\n@two\empty\else _{\box\@ne}\fi}}
\def\t@@left^#1{\@ifnextchar_ {\t@left^{#1}}{\t@left^{#1}_{}}}

\def\two@^#1_#2{\def\n@one{#1}\def\n@two{#2}\mathrel{\setw@dth{#1}{#2}
\mathop{\vcenter{\hbox to \w@dth{\rightarrowfill}\kern-1.7ex
                 \hbox to \w@dth{\rightarrowfill}}%
       }\limits
\ifx\n@one\empty\else ^{\box\z@}\fi
\ifx\n@two\empty\else _{\box\@ne}\fi}}
\def\tw@@^#1{\@ifnextchar_ {\two@^{#1}}{\two@^{#1}_{}}}

\def\tofr@^#1_#2{\def\n@one{#1}\def\n@two{#2}\mathrel{\setw@dth{#1}{#2}
\mathop{\vcenter{\hbox to \w@dth{\rightarrowfill}\kern-1.7ex
                 \hbox to \w@dth{\leftarrowfill}}%
       }\limits
\ifx\n@one\empty\else ^{\box\z@}\fi
\ifx\n@two\empty\else _{\box\@ne}\fi}}
\def\t@fr@^#1{\@ifnextchar_ {\tofr@^{#1}}{\tofr@^{#1}_{}}}

\newdimen\W@dth
\def\setW@dth#1#2{\setbox\z@\hbox{$#1$}\W@dth=\wd\z@
\setbox\@ne\hbox{$#2$}\ifnum\W@dth<\wd\@ne \W@dth=\wd\@ne \fi
\advance\W@dth by 1.2em}

\def\T@^#1_#2{\allowbreak\def\N@one{#1}\def\N@two{#2}\mathrel
{\setW@dth{#1}{#2}
\mathop{\hbox to \W@dth{\rightarrowfill}}\limits
\ifx\N@one\empty\else ^{\box\z@}\fi
\ifx\N@two\empty\else _{\box\@ne}\fi}}
\def\T@@^#1{\@ifnextchar_ {\T@^{#1}}{\T@^{#1}_{}}}

\def\T@left^#1_#2{\def\N@one{#1}\def\N@two{#2}\mathrel{\setW@dth{#1}{#2}
\mathop{\hbox to \W@dth{\leftarrowfill}}\limits
\ifx\N@one\empty\else ^{\box\z@}\fi
\ifx\N@two\empty\else _{\box\@ne}\fi}}
\def\T@@left^#1{\@ifnextchar_ {\T@left^{#1}}{\T@left^{#1}_{}}}

\def\Tofr@^#1_#2{\def\N@one{#1}\def\N@two{#2}\mathrel{\setW@dth{#1}{#2}
\mathop{\vcenter{\hbox to \W@dth{\rightarrowfill}\kern-1.7ex
                 \hbox to \W@dth{\leftarrowfill}}%
       }\limits
\ifx\N@one\empty\else ^{\box\z@}\fi
\ifx\N@two\empty\else _{\box\@ne}\fi}}
\def\T@fr@^#1{\@ifnextchar_ {\Tofr@^{#1}}{\Tofr@^{#1}_{}}}

\def\Two@^#1_#2{\def\N@one{#1}\def\N@two{#2}\mathrel{\setW@dth{#1}{#2}
\mathop{\vcenter{\hbox to \W@dth{\rightarrowfill}\kern-1.7ex
                 \hbox to \W@dth{\rightarrowfill}}%
       }\limits
\ifx\N@one\empty\else ^{\box\z@}\fi
\ifx\N@two\empty\else _{\box\@ne}\fi}}
\def\Tw@@^#1{\@ifnextchar_ {\Two@^{#1}}{\Two@^{#1}_{}}}

\def\to{\@ifnextchar^ {\t@@}{\t@@^{}}}
\def\from{\@ifnextchar^ {\t@@left}{\t@@left^{}}}
\def\two{\@ifnextchar^ {\tw@@}{\tw@@^{}}}
\def\tofro{\@ifnextchar^ {\t@fr@}{\t@fr@^{}}}
\def\To{\@ifnextchar^ {\T@@}{\T@@^{}}}
\def\From{\@ifnextchar^ {\T@@left}{\T@@left^{}}}
\def\Two{\@ifnextchar^ {\Tw@@}{\Tw@@^{}}}
\def\Tofro{\@ifnextchar^ {\T@fr@}{\T@fr@^{}}}

\makeatother

\newcommand{\self}
{   \setlength{\unitlength}{1ex}             
    \begin{picture}(1.5,1.5)(-0.2,-0.1)      
      \put(0.6,0.6){\oval(1.1,0.9)[r]}       
      \put(0.6,0.6){\oval(1.1,0.9)[tl]}      
      \put(-0.35,0.30){\tiny v}              
    \end{picture}                            
}

\def\cahiers#1{Cah\-iers de Topo\-lo\-gie et G\'e\-o\-m\'e\-trie
  Dif\-f\'e\-ren\-ti\-elle\ifnum#1>25 { Cat\-\'e\-go\-rique}\fi, {\bf #1}}

\input diagrams
\newcommand{\hilb}{\mbox{${\bf Hilb}$}}
\newcommand{\cO}{\mbox{$\Omega$}}

\newcommand{\cR}{{\cal R}}

\newcommand{\rel}{{\bf Rel}}
\newcommand{\cS}{{\cal S}}
\newcommand{\cT}[1]{{\cal T}(#1)}
\newcommand{\cD}[1]{{\cal D}(#1)}
\newcommand{\cDp}[1]{{\cal D'}(#1)}
\newcommand{\cE}[1]{{\cal E}(#1)}
\newcommand{\cEp}[1]{{\cal E'}(#1)}

\newcommand{\rarr}{\rightarrow}

\newcommand{\pmon}{{\bf PInj}}
\newcommand{\drel}{{\bf DRel}}
\newcommand{\prel}{{\bf PRel}}
\newcommand{\mrel}{{\bf MRel}}
\newcommand{\cC}{\mbox{${\cal C}$}}

\newcommand{\cLloc}[1]{{\sf L}_{loc}(#1)}
\newcommand{\cH}{\mbox{${\cal H}$}}
\newcommand{\cK}{\mbox{${\cal K}$}}
\newcommand{\reals}{\mbox{${\sf R}$}}
\newcommand{\nuc}[2]{\mbox{${\cal N}(#1,#2)$}}
\newcommand{\meas}{\mbox{${\bf Meas}$}}
\newcommand{\cpd}{\mbox{${\bf Stoch}$}}
\begin{document}
\bibliographystyle{plain}
\title{Nuclear and Trace Ideals in Tensored $*$-Categories}

\author{\begin{tabular}[t]{c}
        Samson Abramsky\thanks{
Research supported in part by EPSRC
}\\
        {\small Department of Computer Science}\\
        {\small University of Edinburgh}\\
        {\small Edinburgh, Scotland}\\
    \end{tabular}
    $\qquad$
        \begin{tabular}[t]{c}
        Richard Blute\thanks{
Research supported in part by NSERC.
}\\
        {\small Department of Mathematics and Statistics}\\
        {\small University of Ottawa}\\
        {\small Ottawa, Ontario, Canada}\\
    \end{tabular}
     \\ 
        \begin{tabular}[t]{c}
        \rule{0mm}{10mm}\\
        Prakash Panangaden\thanks{
Research supported in part by NSERC.  Much of this work was done while on 
leave at BRICS, Aarhus University. 
}\\     
        {\small Department of Computer Science}\\
        {\small McGill University}\\
        {\small Montreal, Quebec, Canada}\\
    \end{tabular}
    }
\date{June 7th, 1998}
\maketitle

\begin{center}
{\em Presented to Mike Barr on the occasion of his 60th birthday.}
\end{center}
\begin{abstract}
We generalize the notion of nuclear maps from functional analysis by
defining nuclear ideals in tensored $*$-categories.  The motivation for
this study came from attempts to generalize the structure of the category
of relations to handle what might be called ``probabilistic relations''.
The compact closed structure associated with the category of relations 
does not generalize directly, instead one obtains nuclear ideals.

Most tensored $*$-categories have a large class of morphisms
which behave as if they were part of a compact closed category, i.e. they
allow one to transfer variables between the domain and the codomain.  We
introduce the notion of {\em nuclear ideals} to analyze these classes of
morphisms.  In compact closed tensored $*$-categories, all morphisms
are nuclear, and in the tensored $*$-category of Hilbert spaces, 
the nuclear morphisms are the {\it Hilbert-Schmidt maps.}

We also introduce two new examples of tensored $*$-categories, in which
integration plays the role of composition. In the first, morphisms are a 
special class of {\em distributions}, which we call {\it tame distributions}. 
We also introduce a category of {\em probabilistic relations}.

Finally, we extend the recent work of Joyal, Street and Verity 
on {\it traced monoidal categories} to this setting
by introducing the notion of a {\it trace ideal}. 
We establish a close correspondence between nuclear ideals and trace 
ideals in a tensored $*$-category, suggested by the correspondence 
between Hilbert-Schmidt operators and 
trace operators on a Hilbert space. 

\end{abstract}

\pagebreak

\section{Introduction}
This paper develops a new categorical structure, called a {\em nuclear
ideal}, which comes from two independent, seemingly unrelated,
developments.  These are Grothendieck's concept of nuclearity in
functional analysis, see for example~\cite{Treves67},
and the usual notion of binary relations.  The original motivation for
this investigation was the need to generalize ordinary binary relations to
probabilistic relations with an eye towards certain applications in
computer science.  However a satisfactory notion of what this
generalization should be comes from the concept of nuclearity in functional
analysis.  This paper presents the new concept and gives several nontrivial
examples of nuclear ideals.

Relations form a basic and ubiquitous mathematical structure.  There has
been much activity in formulating what relations are ``abstractly'', so
that one can generalize the concept to new situations.  Typical examples of
such formulations are the concept of cartesian bicategories~\cite{Carboni87}
and allegories~\cite{Freyd90}.  One of the key aspects of the category
\rel{} is the fact that one has ``transfer of variables''
i.e. one can use the closed structure and the involution to move variables
from ``input'' to ``output''.  Intuitively speaking, this reflects the idea
that the source and target of a binary relation are a matter of convention
and a binary relation is an inherently symmetric object.  In many
situations that otherwise resemble relations, one finds that the closed
structure does not exist and hence one loses the ability to transfer variables.
A typical analogue of binary relations are the ``probabilistic'' 
binary relations, described at length later in
the paper.  Even in the absence of detailed definitions it ought to be
clear that one cannot (indeed {\em should not}) rearrange the inputs and
outputs of a probabilistic relation because there may be dependencies 
present among different inputs.  What remains then in lieu of closed
structure?  We claim that it is precisely the nuclear ideals of the present
paper.

In these settings, there appears to be a tension between having 
identities and having compact closed structure. If one looks only at the 
nuclear ideal, one has a compact closed ``category'' without identities. 
On the other hand, the ambient category lacks closed structure. Others have
observed that there are ``categories without identities'', and given a wide
range of examples and applications \cite{Ageron96,Selinger97}. However the
interplay between the ideal and the ambient category is the point
of the present work, not just the lack of identities. 

Another motivation for this work comes from considering Hilbert spaces.
The tensored $*$-category of Hilbert spaces and bounded linear maps 
(hereafter denoted ${\bf Hilb}$)
shares much of the same structure as \rel{}.  One of the goals of
this paper is to measure the extent of this correspondence.  Like the
category of relations, ${\bf Hilb}$\ has a tensor product and a
tensor-preserving involution, which is the identity on objects.  In the
case of ${\bf Hilb}$, it is given by the adjoint operation.  However, the
category of Hilbert spaces lacks the closed structure of \rel{}.  The
structure of ${\bf Hilb}$\ has been axiomatized as the notion of a {\it
tensored $*$-category} \cite{Ghez85,Doplicher89}.  (In fact, it
is a {\it tensored $C^*$-category}, but we will not consider its normed
structure here.)

In this paper, we argue that a tensored $*$-category should be thought of as a 
category of (generalized) relations. The category of relations is compact 
closed, and this property is frequently taken to be fundamental in axiomatizing
relational categories \cite{Abramsky94,Carboni87}.
However, the categories of relations which we consider are not compact closed, 
but rather contain a large class of morphisms, in fact an ideal, which 
has the basic structure of a compact closed category.  
To axiomatize this notion, we introduce the new notions of
{\it nuclear ideal} and {\it nuclear morphism.} This idea is based
on the definition of a nuclear morphism between Banach spaces, due to
Grothendieck~\cite{Grothendieck55}, which was subsequently axiomatized by
Higgs and Rowe~\cite{Higgs89}.  
The concept of nuclearity in analysis can be viewed as describing when one
can think of linear maps as matrices.  Of course, in the finite-dimensional
case one can always do this and it will be the case that all maps between
finite-dimensional vector spaces are nuclear.  
The Higgs-Rowe theory applies only to
autonomous (symmetric monoidal closed) 
categories, while our definition applies to the somewhat
different setting of tensored $*$-categories.  In the case of a compact closed
$*$-category, all morphisms are nuclear, while in ${\bf Hilb}$ with its
usual tensored $*$-structure, the
nuclear morphisms are precisely the Hilbert-Schmidt maps 
\cite{Kadison83}. Note that since we are only considering ${\bf Hilb}$
with the $L_2$ tensored $*$-structure, the notion of nuclear map we obtain is 
different from Grothendieck's notion arising from the category of Banach
spaces (with, of course, the $L_1$ tensor product).

A further goal of this paper is to introduce two new examples of tensored
$*$-categories, in which integration plays the role of composition.  The
first such category is a category of {\it generalized functions} or {\it
distributions} \cite{Algwaiz92,Treves67}.  Since a discrete relation on
$X\times Y$ can be viewed as a function $f:X\times Y\rarr \{0,1\}$, it
seems reasonable to model a ``smeared out'' relation as a continuous
function $f:U\times V\rarr {\sf R}$, where $U$\ and $V$\ are open subsets
of Euclidean space.  However, the identity for such a category would be the
{\it Dirac Delta} which is not a function, but a distribution.  
We choose a particular class of distributions, the {\it
tame distributions}, which are sufficiently functional to allow
composition. We then present a nuclear ideal for this category. It will consist
of the tame distributions with functional kernel. 

To build a category of {\it probabilistic relations}, one would like a
category where the objects are probability spaces, and a morphism is a
measure on the product space.  The structure we eventually arrive at is the
notion of {\em conditional probability distribution}, described in section
9.  Categories of conditional probability distributions have previously
been studied by Giry~\cite{Giry80} and Wendt~\cite{Wendt93,Wendt94}.  Our
formulation differs from theirs in that in our category, objects 
are equipped 
with measures and morphisms are measures on the product space satisfying an
absolute continuity property. To each morphism, we are then able to associate
a pair of conditional probability distributions.
Again, in this case the nuclear ideal will consist 
of measures having a functional kernel.

We also extend the recent work of Joyal, Street and Verity 
on {\it traced monoidal categories} \cite{Joyal96} to the present setting
by introducing the notion of a {\it trace ideal}. For a given symmetric 
monoidal category, it is not generally the case that
arbitrary endomorphisms can be assigned a trace. However, one can often
find ideals on which a trace can be defined satisfying equations analogous
to those of Joyal, Street and Verity. Our abstract definition is suggested
by the usual trace construction in the category of Hilbert spaces, where
there is a well-established relationship between maps in the trace class 
and Hilbert-Schmidt maps. In this case, we obtain the usual notion of trace
of a bounded linear operator in the trace class.

\bigskip

\noindent {\bf Acknowledgments}- The authors would like to thank Mayer Alvo,
John Baez, Michael Barr, Robin Cockett, Thomas Ehrhard, Jean-Yves Girard, 
Martin Hyland, Vojkan Jaksic, William Lawvere, Robert Seely and 
Michael Wendt for helpful discussions. We also received a number of extremely
helpful comments from the anonymous referees.

\section{Categorical Preliminaries}

We assume the reader is familiar with the notion of a {\it
symmetric monoidal\footnote{We observe that the 
notation $c\colon A\ox B\rarr B\ox A$ is used for the symmetry, and $I$ is 
used for the tensor unit.} 
category}. A suitable reference is ~\cite{MacLane71}.
We now review some of the different closed structures such a category could
have.

\begin{defn}{\em A symmetric monoidal category is {\it closed} or {\it
autonomous} if, for all objects $A$\ and $B$, there is an object $A\oto B$
and an adjointness relation:

\[ Hom(A\ox B, C)\cong Hom(B, A\oto C)\]

The unit and counit of this adjunction are the familiar morphisms:

\[ ev\colon A\ox (A\oto B)\rarr B \ \ \ \ coev\colon A\rarr B\oto(A\ox B)\]

}\end{defn}

Examples of autonomous categories include the category of vector spaces and 
the category of relations. We obtain a pair of autonomous categories by 
considering Banach spaces. We can either consider ${\bf Ban}_{\infty}$, the
category of Banach spaces and bounded linear maps, or we can consider
the category ${\bf Ban}_{1}$, of Banach spaces and maps of norm less than
or equal to 1. In either case, the internal Hom is the Banach 
space of all bounded linear maps, and the tensor product is the completed 
projective tensor product ~\cite{Treves67}.

\begin{defn}{\em
A {\em compact closed category} is a symmetric monoidal category such
that for
each object $A$ there exists a dual object $A^{*}$, and canonical
morphisms:

\begin{center}
$\nu\colon I\rightarrow A\ox A^{*}$\\
$\psi\colon A^{*}\ox A\rightarrow I$
\end{center}

\noindent such that the usual adjunction equations hold:

$$
\begin{diagram}
&&I\ox A\\
&\NE^{\cong}&&\SE^{\nu\ox id}\\
A&&&&A\ox A^{*}\ox A\\
&\rdTo(1,2)_{id}&&\ldTo(1,2)_{id\ox \psi}\\
&A &\lTo_{\cong}& A\ox I\\
\end{diagram}
$$

\noindent together with the dual diagram for $A^*$. In the case of
a strict monoidal category, these equations reduce to the usual adjunction
triangles. It is easy to 
see that a compact closed category is indeed closed and that
$A\oto B\cong A^{*}\ox B$.
} 
\end{defn}

Compact categories could also be defined as $*$-autonomous 
categories~\cite{Barr80,Blute93}  with
the additional canonical isomorphism $A^{*}\ox B^{*}\cong (A\ox B)^{*}$.
$*$-Autonomous categories provide the basic framework for the model theory of 
the multiplicative fragment of linear logic \cite{Girard87}.

We briefly describe the prototypical example,  the {\it 
category of relations}.

\begin{defn}{\em
The {\em category of relations}, 
\rel{}, has sets as objects, a morphism from $X$
to $Y$ will be a relation on $X\times Y$, with 
the usual relational composition.
}\end{defn}

In what follows, $X,Y,Z$ will denote sets, and $x,y,z$ will denote elements.
A binary relation on $X\times Y$ will be denoted $x\cR y$. The identity
relation will be denoted ${\cal ID}$, and is defined as $x{\cal ID}x$, for
all $x\in X$. Given a relation $\cR\colon X\rarr Y$, we let 
$\overline{\cR}\colon Y\rarr X$
denote the converse relation.

We verify that \rel{} is compact. 
The tensor product $\ox$ is given by taking the products of sets, and 
on morphisms, we have:

\[ \cR\colon X\rarr Y \;\;\;\;\;\;\;\;\; \cS\colon X'\rarr Y' \]

\begin{center} 
$(x,x')\cR\ox\cS (y,y')$ if and only if $x\cR y$ and $x'\cS y'$
\end{center}

The unit for the tensor is given by any one point set.
We define the functor $(\;)^{*}\colon \rel\rarr\rel{}$ by:

\[ X^{*}=X \;\;\;\;\; \cR^{*}=\overline{\cR} \]

The relation $\nu\colon I\rightarrow X\ox X^{*}$ is given by $*\nu (x,x)$
for all $x\in X$ and similarly for $\psi$.

\section{The Tensored $*$-Category of Hilbert Spaces}

Our notation for this section will be as follows.  
We will use brackets of the form $\langle -,- \rangle$\ to 
denote the inner product, which will be linear in the first variable.  
The associated norm will be denoted $||-||$.  If $\alpha$ is an element of
the base field, then $\overline{\alpha}$ will denote its conjugate. If $H$
is a Hilbert space, then $\overline{H}$ will denote the conjugate space.
An orthonormal basis will be denoted $\{e_{i}\}_{i\in I}$.
A suitable reference for basic Hilbert space theory is \cite{Kadison83}.

Let \hilb \ denote the category of Hilbert 
spaces and bounded linear maps, where ``bounded'' always 
means bounded in the norm associated to the inner product. 
We now discuss the structure of this category 
which is relevant to this paper.  The 
first structure we need is the {\it adjoint function} \cite{Kadison83}.

\begin{defn}{\em 
Let \cH \ and \cK \ 
be Hilbert spaces, and $f\colon \cH\rarr \cK$\ a bounded linear map.
Then the {\it adjoint} of $f$, denoted $f^*$, is defined to be the unique
bounded linear map
$f^{*}\colon\cK\rarr \cH$\ such that, for all $a\in \cH, \ \ b\in \cK$, we
have:

\[ \langle a,f^{*}(b)\rangle =\langle f(a),b\rangle\]

}
\end{defn}

\begin{lem}
The adjoint construction satisfies the following properties:

\begin{itemize}
\item $(id_{\cH})^{*}=id_{\cH}$
\item $(fg)^{*}=g^{*}f^{*}$
\item $f^{**}=f$
\item $(f\ox g)^{*}=f^{*}\ox g^{*}$ (The tensor product will be discussed 
below.)
\end{itemize}
\end{lem}

These conditions tell us that the adjoint operation provides a 
contravariant, tensor-preserving, 
involutive functor on \hilb\ which is the identity on objects.
Given such a functor, it is clear that the category \hilb\ is 
much closer in its categorical
structure to the category of relations than to the category of 
Banach spaces. 

\subsection{Hilbert-Schmidt Maps}

We now discuss a crucial class of bounded linear maps, called the {\it
Hilbert-Schmidt maps}.  The material in this section can be found in
\cite{Kadison83}.

\begin{defn}{\em
If $f\colon\cH\rarr \cK$ is a bounded linear map, we call $f$\
a {\it Hilbert-Schmidt map} if the sum

\[ \sum_{i\in I} ||f(e_{i})||^{2}\]

\noindent is finite for an orthonormal basis $\{e_{i}\}_{i\in I}$.  
The sum is independent of basis chosen.
}\end{defn}

It is straightforward to see that:

\begin{lem}
If  $f\colon\cH\rarr \cK$ \ is a Hilbert-Schmidt map and 
$g\colon\cH_{1}\rarr \cH$, 
$g'\colon \cK\rarr \cK_{1}$\ are arbitrary bounded linear maps, then 
$g'f$\ and $fg$\ are Hilbert-Schmidt.
\end{lem}

Thus the Hilbert-Schmidt operators on a space form a 2-sided
ideal in the set of all bounded linear operators.  
A proof of the following theorem may be found in \cite{Kadison83}.  

\begin{thm}
Let ${\bf HSO}(\cH,\cK)$\ denote the set of Hilbert-Schmidt maps from \cH\ to
\cK.  Then ${\bf HSO}(\cH,\cK)$\ is a Hilbert space with:

\[\langle f,g\rangle = 
\sum_{i\in I, j\in J} \langle f(e_{i}),e'_{j}\rangle \langle
e'_{j},g(e_{i})\rangle \]
\end{thm}

Here, $\{e_{i}\}_{i\in I}$ is an orthonormal basis for \cH \ and
$\{e_{j}'\}_{j\in J}$ is an orthonormal basis for \cK.

\subsection{The Tensor Product}

It is standard to construct the tensor product of Hilbert spaces $\cH\ox\cK$\
as the completion of the algebraic tensor product with respect to the inner 
product:

\[ \langle x_{1}\ox y_{1}, x_{2}\ox y_{2}\rangle
=\langle x_{1},x_{2}\rangle \langle y_{1},y_{2}\rangle \]

One then completes with respect to the $L_2$ norm to obtain a Hilbert space.
(Note that it is also possible 
to give an equivalent presentation that emphasizes the universal 
mapping property of the tensor. This involves the notion of a {\it
weak Hilbert-Schmidt mapping}. This is explained in \cite{Kadison83}, page
132.)

\begin{rem}
We wish to emphasize that, in this paper, we will only be considering the $L_2$
tensor product. Furthermore, the category ${\bf Hilb}$ will always be
the category of Hilbert spaces with bounded linear maps, equipped 
with this tensored $*$-structure. 

Thus, our notion of nuclearity will not coincide with the notion obtained 
by viewing Hilbert spaces as Banach spaces and applying Grothendieck's 
definition, which of course uses the $L_1$ tensor.
\end{rem}

For us, the most important property of the Hilbert tensor product is 
its relation to Hilbert-Schmidt maps. This is given by the following 
theorem \cite{Kadison83}, p.142:

\begin{thm}\label{nuctheorem}
We define a linear mapping $U\colon\overline{\cH}\ox\cK\rarr {\bf HSO}
(\cH,\cK)$\
by \mbox{$U(x\ox y)(u)=\langle x,u\rangle y$}, where $x\ox y\in\cH\ox\cK$.  
Then $U$\ is a unitary transformation of 
$\overline{\cH}\ox\cK$ onto ${\bf HSO}(\cH,\cK)$. In particular, we
note that the morphism $U$ is a linear bijection.
\end{thm}

\section{Tensored $*$-categories}

The category \hilb, of Hilbert spaces and bounded linear maps, 
shares many of the properties of a compact closed 
category, except for the closed structure.  \hilb\
is in fact an example of a {\it tensored $*$-category}. 
We now develop this theory.

\begin{defn}
{\em 
A category $\cal C$\ is a {\em $*$-category} if it is equipped with a functor
$(-)^*\colon{\cal C}^{op}\rarr{\cal C}$, which 
is strictly involutive and the identity
on objects.  (Note that the strict involution may be replaced with a coherent 
involution, but we will not require that level of generality.) 
A $*$-category is {\em tensored} if it is symmetric monoidal, 
$(f\ox g)^{*}=f^{*}\ox g^{*}$, and there is a 
covariant {\it conjugate functor}, 
$\overline{(\,\,)}\colon{\cal C}\rarr {\cal C}$, which commutes 
with the $*$-functor and has natural isomorphisms: 

\begin{itemize}
\item $\overline{\overline{A}}\cong A$ (We will generally take this to be an 
equality.)
\item $\overline{A\ox B}\cong\overline{A}\ox\overline{B}$ (We 
will generally take this to be an equality.)
\item $\overline{I}\cong I$. 
\end{itemize}

\noindent satisfying the usual monoidal equations, and the following equation.
Suppose that $f\colon I\rarr I$. 

$$\begin{diagram}
I&\rTo^{f^*}&I\\
\dTo^\cong&&\uTo_\cong\\
\overline{I}&\rTo_{\overline{f}}&\overline{I}\\
\end{diagram}$$

In all of our examples except for those involving 
complex Hilbert spaces, conjugation will simply be taken to be the identity.
In this case, the previous diagram implies 
that if $f\colon I\rarr I$, then $f^*=f$.
}\end{defn}

The notion of a tensored $*$-category is the first step towards defining a 
tensored {\it $C^*$-category}, or multiobject $C^*$-algebra 
\cite{Ghez85,Doplicher89}. This theory has been developed 
quite extensively in the previously cited references. Among the results
established is a representation theorem stating that such categories have 
faithful structure-preserving embeddings in ${\bf Hilb}$. This should be 
thought of as a multiobject version of the Gelfand-Naimark-Segal theorem.

\bigskip

\noindent{\bf Examples of tensored $*$-categories}

\begin{itemize}
\item \rel{} 
\item \hilb\ 
\item $\hilb_{fd}$, the category of finite-dimensional Hilbert spaces.
\item ${\bf URep}(G)$, the category of unitary representations of a 
compact group $G$.
\item ${\bf URep}_{fd}(G)$, the category of finite-dimensional 
unitary representations of a 
compact group $G$.
\end{itemize}

Further examples can be found in \cite{Ghez85,Doplicher89}.
Note that examples 2 and 4 are tensored $*$-categories which are not closed. 
We will present other examples of tensored $*$-categories which are not 
closed.

Even though tensored $*$-categories are not compact closed, they share much
of the same structure.  One of the goals of this paper is to introduce
a structure for measuring the extent to which such a category is closed.

\section{Nuclearity}

One of the characteristic features of compact closed categories is 
the ability to distribute the dual functor across the tensor product. This 
is represented by the equation $(A\ox B)^{\circ}\cong 
A^{\circ}\ox B^{\circ}$. ($A^\circ$ denotes the dual object. We 
temporarily adopt this notation to avoid confusion with the $*$-functor
we will be discussing later. In the context of tensored $*$-categories, 
one should keep in mind the equation $A^{\circ}=\overline{A^*}$.)
This allows one to arbitrarily repartition
the morphism or ``interface'' in the terminology of interaction categories
~\cite{Abramsky94}. The categories we will encounter typically allow such 
repartitioning for some maps, but do not meet all the 
requirements of being a compact closed category. 

We now introduce the related notion of {\it nuclearity} in a symmetric 
monoidal closed category, 
due to Rowe~\cite{Rowe88}, and subsequently studied by Rowe and Higgs 
\cite{Higgs89}. The idea is suggested by Grothendieck's work on topological 
tensor products and nuclear spaces
\cite{Grothendieck55}. Grothendieck defined a continuous linear map
$f\colon A\rarr B$ between Banach spaces to be {\it nuclear} if it can be written
as $f(a)=\Sigma f_i(a)b_{i}$ where $\Sigma f_{i}\ox b_{i}$ is an element
of the completed projective tensor product $A^{\circ}\ox B$. We begin by noting 
that in any symmetric monoidal closed category, there is a morphism of the 
form:

\[ \varphi\colon B\ox A^{\circ}\rarr A\oto B\]

Here, $A^{\circ}=A\oto I$, where $I$\ is the unit for the tensor.

\noindent This is constructed as the transpose of the composite:

\[ B\ox A^{\circ}\ox A\stackrel{id\ox\psi}{\longrightarrow} B\ox I
\stackrel{\cong}{\longrightarrow} B\]

\begin{defn}{\em Let $\cC$ be a symmetric monoidal closed 
category. Let $\varphi$ denote
the canonical morphism \mbox{$\varphi\colon B\ox A^{\circ}\rarr A\oto B$.} If 
$f\colon A\rarr B$ in
$\cC$, then let $n(f)\colon 
I\rarr A\oto B$ be the name of $f$. We say that $f$ 
is {\em nuclear} if there exists $p(f)\colon I\rarr B\ox A^{\circ}$\ 
such that the following diagram commutes:

$$
\begin{diagram}
I &&\rTo^{p(f)}&&B\ox A^{\circ}\\
&\SE_{n(f)}&&\SW_{\varphi}\\
&&A\oto B
\end{diagram}$$

\bigskip

\noindent We will refer to $p(f)$ 
as a {\em pseudoname} for $f$. (We should point
out that there are some cases in which a pseudoname is not unique.)
We say that an object of $\cC$ is {\it nuclear} if its identity map is nuclear.
  
}\end{defn}

\begin{lem}
Suppose that 
$f\colon A\rarr B$\ and $g\colon C\rarr D$\ are nuclear, then so are:
\begin{itemize}
\item $f^{\circ}\colon B^{\circ}\rarr A^{\circ}$
\item $f'f\colon A\rarr E$ for any morphism $f'\colon B\rarr E$
\item $fh\colon F\rarr B$ for any morphism $h\colon F\rarr A$
\item $f\ox g\colon A\ox C\rarr B\ox D$
\end{itemize}
\end{lem}

All of the above can be obtained by straightforward diagram chasing. For 
example, in the third item, one can choose $p(fh)=p(f);(id \ox h^{\circ})$.
It is not in general the case that if $f,g$ are nuclear, then so is $f\oto g$.
However, if $\cal C$ is $*$-autonomous with unit as dualizing object, then
$f\oto g$ will also be nuclear, \cite{Higgs89} p. 70.

In a compact closed category, the map $\varphi$\ is an isomorphism, and thus 
every 
map is nuclear. Furthermore, we can see the following:

\begin{thm}[\cite{Higgs89}, Thm. 2.5]
For an arbitrary object $A$ in $\cal C$, a symmetric monoidal closed 
category, the following are equivalent.

\begin{itemize}
\item $A$ is nuclear
\item The morphism $\varphi\colon A\ox A^\circ\rarr A\oto A$ is an isomorphism.
\item The morphism $\varphi\colon B\ox A^\circ\rarr A\oto B$ is an isomorphism,
for arbitrary objects $B$.
\end{itemize}
\end{thm}

\begin{thm}
For any symmetric monoidal closed category, the full subcategory of nuclear
objects is compact-closed.
\end{thm}

\begin{prf}{}

Suppose that $A$\ is a nuclear object. Then, choosing a pseudoname for the 
identity gives a morphism of the form $I\rightarrow A\ox A^{\circ}$. It only 
remains to show that the adjunction triangles commute. We will consider one
of the two adjunction triangles.

$$
\begin{diagram}
A\cong I\ox A & \rTo^{p(id)\ox id} & A\ox A^{\circ}\ox A\\
\dTo^{n(id)\ox id} & \SW^{\varphi\ox id} & \dTo_{id\ox ev}\\
(A\oto A)\ox A&\rTo_{ev}&A\ox I\cong A\\
\end{diagram}$$

It is standard that the lower leg of the diagram is the identity. The upper
leg of the diagram corresponds to the adjunction triangle. The upper triangle
in the above square is the definition of pseudoname. The lower triangle is a
straightforward exercise.

\end{prf}

In ${\bf Ban}_{\infty}$ or ${\bf Ban}_{1}$, we recover Grothendieck's 
original definition of nuclearity. The nuclear objects 
are the finite-dimensional Banach spaces.
In the category of vector spaces, a morphism is nuclear if and only if
its image is finite-dimensional. Again, a vector 
space is nuclear if and only if it is finite-dimensional.
In \cite{Higgs89}, the authors explore the notion of nuclearity in the category
of complete join semilattices ${\bf CJSL}$. 
It is well known that this is a symmetric monoidal closed
category, in fact $*$-autonomous \cite{Joyal84}. The authors completely 
characterize nuclearity in this case (This 
result is closely related to Raney's notion of a {\it tight} morphism
\cite{Raney60}.):

\begin{thm}(Higgs, Rowe)
A morphism $f\colon A\rarr B$ in ${\bf CJSL}$ is 
nuclear if and only if there exists $g\colon B\rarr A$
such that for all $a\in A, f(a)=sup\{b|a\not\leq g(b)\}$. An object is
nuclear if and only if it is completely distributive.
\end{thm}

\begin{rem}
Following recent work of Joyal, Street and Verity \cite{Joyal96} on 
{\it traced monoidal categories}, one can now observe that, 
in a symmetric monoidal closed category, it is possible to 
define a trace on the nuclear morphisms as follows, under the assumption 
that pseudonames are unique. If $f\colon A\rarr A$ is
nuclear, then $tr(f)\colon I\rarr I$ is given by:

\[ tr(f)=p(f);ev\colon I\rarr A\ox A^\circ\rarr I\]

\noindent where $ev\colon A\ox A^\circ\rarr I$ is the usual 
evaluation map. Then given $h\colon A\rarr B$ a nuclear map, and 
$g\colon B\rarr A$ arbitrary, one can verify the usual trace equation
$tr(gh)=tr(hg)$. This is seen by the following diagram:

$$
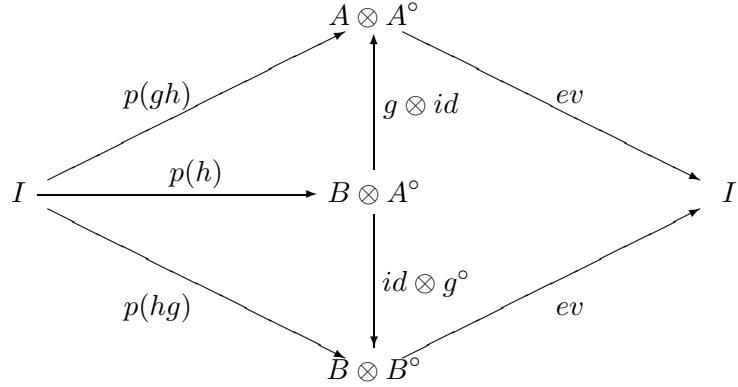
\begin{diagram}
&&&&A\ox A^\circ\\
&&&\ENE^{p(gh)}&\uTo_{g\ox id}&\ESE^{ev}\\
I&&\rTo^{p(h)}&&B\ox A^\circ&&&&I\\
&\ESE_{p(hg)}&&&\dTo_{id\ox g^\circ}&&&\ENE_{ev}\\
&&&&B\ox B^\circ\\
\end{diagram}$$

\bigskip

The righthand diamond is the usual (di)naturality of evaluation. The two 
triangles on the left are the equations for $p(gh)$ and $p(hg)$. 

\end{rem}

While this theory is satisfactory when considering symmetric monoidal closed 
categories,
there are nonclosed categories which exhibit similar structure. For example, 
the category of Hilbert spaces is not closed, but the 
class of Hilbert-Schmidt maps seem to have
something like a nuclearity property. We will soon exhibit other
such categories. One of the goals of this paper is to extend 
the above notions to a larger class of categories, specifically to 
$*$-categories. We now introduce a new notion, that of a {\it 
nuclear ideal}.

\begin{defn}
{\em

Let $\cal{C}$ be a tensored $*$-category. A {\em nuclear ideal} for 
${\cal{C}}$ consists of the following structure: 

\begin{itemize}
\item For all objects $A, \, B\in \cC$, a subset 
${\cal{N}} (A,B)\subseteq Hom(A,B)$. We will refer to the union of these 
subsets as ${\cal{N}} (\cal{C})$\ or $\cal{N}$. We will refer to the
elements of ${\cal{N}}$\ as {\em nuclear maps}. The class $\cal{N}$\ 
must be closed under composition with arbitrary 
\cC-morphisms, closed under $\ox$, closed under $(\,)^{*}$, and the conjugate 
functor.
\item A bijection $\theta\colon {\cal{N}}(A,B) \rarr Hom(I,\overline{A}\ox B)$.
If $f\colon A\rarr B$ is a nuclear morphism, note that we can 
use the bijection $\theta$ and the $*$-functor to construct morphisms of the 
form:

\begin{enumerate}
\item $\theta(f)\colon I\rarr \overline{A}\ox B$
\item $\theta(f)^*\colon\overline{A}\ox B\rarr I$ 
\item $\theta(f^*)\colon I\rarr \overline{B}\ox A$
\item $\theta(f^*)^*\colon \overline{B}\ox A\rarr I$ 
\end{enumerate}

\noindent 
We shall frequently refer to these morphisms as {\it transposes} 
of $f$. It will always be clear from the
context which transpose is being considered. The bijection $\theta$ must also 
satisfy the following properties:

\begin{enumerate}

\item {\bf Preservation of tensored $*$-structure} \
The bijection $\theta$ must preserve all of the tensored $*$-structure. 
In other words,

\smallskip

\begin{enumerate}
\item If $f:A\rarr B$ and $g:C\rarr D$ are nuclear, then $\theta(f\ox g)=
\theta(f)\ox \theta(g)$. More precisely, the map $\theta(f\ox g)\colon
I\rarr\overline{A\ox C}\ox B\ox D$ is given by the composite:

\[I\cong I\ox I\rarr \overline{A}\ox B\ox\overline{C}\ox D\cong
\overline{A}\ox\overline{C}\ox B\ox D=\overline{A\ox C}\ox B\ox D\]

Furthermore, the transposes of a map of the form $f\colon I\rarr A$ are 
given by composition with the evident isomorphism.

\smallskip

\item $\theta(\overline{f})=\theta(f^*)=\overline{\theta(f)}$. 
Again, more precisely, we would say:

\[\theta(\overline{f})=c\circ\theta(f^*)=\overline{\theta(f)}\circ\iota\]

\noindent where $c$ is the symmetry and $\iota$ is the isomorphism
$\iota\colon I\rarr \overline{I}$.

\end{enumerate}

\smallskip

\item {\bf Naturality} For any 
$f\colon A\rarr C$\ and $g\colon B\rarr D$, the following 
diagram commutes:

$$\begin{diagram}\label{nucdiagram}
{\cal{N}} (A,B)&\rArr^{\theta}&Hom(I,\overline{A}\ox B)\\
\dArr^{{\cal{N}} (f^*,g)}&&\dArr_{Hom(I, \overline{f}\ox g)} \\
{\cal{N}} (C,D)&\rArr^{\theta}&Hom(I,\overline{C}\ox D)\\
\end{diagram}$$

\noindent Note that since the class of nuclear morphisms is closed 
under composition with arbitrary \cC-morphisms, the function 
${\cal{N}} (f^*,g)$ \ is well defined.

\smallskip

\item {\bf Compactness} Let $f\colon A\rarr B$\ and $g\colon 
B\rarr C$\ be nuclear.

Then the following should commute.

$$\begin{diagram}
A&\rArr^{\cong}&I\ox A&&\rArr^{\theta(g)\ox id_A}&&\overline{B}\ox C\ox A\\
\dArr_{gf}&&&&&&\dArr^{c}\\
C&\lArr_{\cong}&C\ox I&&\lArr_{id_C\ox\theta(f^*)^*}&&C\ox\overline{B}\ox A\\
\end{diagram}$$

\end{enumerate}

\end{itemize}

This completes the definition of nuclear ideal. 
In the case where $A$ is a nuclear object and $f=g=id_A$, then this last
equation reduces to the usual adjunction equation for a compact closed 
category. We will see that it is also related to the ``yanking'' axiom of
\cite{Joyal96}.
}
\end{defn}

\bigskip

Given a category $\cal C$ and a nuclear ideal  $\cal N$, we say
that an object $A$ of $\cal C$\ is $\cal N$-nuclear if we have that 
${\cal{N}}(A,-)=Hom(A,-)$. Note that by the ideal property, this is equivalent
to saying that the identity map for $A$\ is nuclear. Typically, this notion
of nuclear object is capturing the ``finite-dimensional'' subcategory. It
should not be thought of as describing Grothendieck's much richer theory
of nuclear spaces.

Note that we are not claiming that the transposition map is in any way 
unique; different choices of $\theta$ could conceivably
give different nuclear ideal structures. The usual uniqueness 
arguments, see for example \cite{MacLane71} pp. 80-82, do not apply here 
in that we may not transpose the identity
map. Thus it is possible that several distinct nuclear structures may exist
on a given category. We are still pursuing this question. 
However, we know of no such examples.  
In the examples presented in this paper, the choice of the
transpose is obvious and canonical, given the structures 
under consideration.

One of the consequences of the above definition is the ``sliding'' equation of 
Joyal, Street and Verity \cite{Joyal96}:

\begin{lem} Suppose $f\colon A\rarr B$ and $g\colon B\rarr A$ are nuclear.
Then the following diagram commutes for any nuclear ideal:

$$\begin{diagram}
&&\overline{A}\ox B\\
&\NE^{\theta(f)}&&\SE^{\theta(g^*)^*}\\
I&&&&I\\
&\SE_{\theta(g)}&&\NE_{\theta(f^*)^*}\\
&&\overline{B}\ox A\\
\end{diagram}$$
\end{lem}

This equation is a straightforward consequence of the axioms. We will see 
in Section 8 that it corresponds to the familiar trace equation 
$tr(fg)=tr(gf)$.

\begin{thm}
Let $({\cal C},{\cal N})$\ be a nuclear ideal for
which all objects are nuclear, then ${\cal C}$ is a 
compact-closed category.
\end{thm}

\begin{prf}{}

If $A$\ is an object of $\cal C$, then the transpose of the identity will be
a morphism of the form \mbox{$I\rarr A\ox \overline{A}$}. The commutativity
of the adjunction triangles follows from the compactness requirement of the
definition.

\end{prf}

\begin{thm}
The set of Hilbert-Schmidt maps forms a nuclear ideal for \hilb.
\end{thm}

\begin{prf}{}

Let ${\cal H}$\ and ${\cal K}$\ be Hilbert spaces, and let 
${\cal N}({\cal H},{\cal K})$ be the set of all Hilbert-Schmidt maps from 
${\cal H}$\ to ${\cal K}$. It is evident that
$Hom(I,\overline{\cH}\ox\cK)\cong\overline{\cH}\ox\cK$. 
So the morphism $U$, defined in \ref{nuctheorem}, 
will act as a transpose operator. We saw in section 3 that this map 
was a linear bijection. It only remains to check the equations.  These are
 a straightforward consequence of linearity and 
properties of the adjointness operator.

\end{prf}

The nuclear objects in this case are precisely the 
finite-dimensional Hilbert spaces. Thus we recover the familiar compact closed
subcategory. The same program can be carried out for categories of 
representations such as ${\bf URep}(G)$.

\subsection{Partial Injective Functions}

Define a category \pmon\ as follows. Its objects will be sets, and morphisms
will be partial injective functions, that is to say partial functions which are
monomorphic when restricted to the domain. These partial functions 
were used by Danos in his modeling of the geometry of interaction 
\cite{Danos}.

If $f\colon X\rarr Y$ is a morphism,
let $Dom(f)$ be its domain, i.e. $Dom(f)=\{x\in X|f(x) \mbox{ is defined}\}$. 
This category has an evident $*$-structure, and if we choose the cartesian
product of sets as a tensor, then we evidently have a tensored $*$-category.
We now demonstrate that this category has an evident nuclear ideal. Define:

\[{\cal N}(X,Y)=\{f\colon X\rarr Y|Dom(f)\mbox{ has cardinality 0 or 1}\}\]

Then one can see that we have an obvious bijection between
$Hom(I,X\ox Y)$ and ${\cal N}(X,Y)$.

\begin{thm}
The above construction defines a nuclear ideal for \pmon.
\end{thm}

\subsection{Crossed $M$-Sets}

The following is based on Freyd and Yetter's notion of a {\it crossed G-Set}, 
which they use in their work on braided compact closed categories 
\cite{Freyd89}. In this paper, we will only consider 
a commutative monoid, which gives a symmetric monoidal category.
We hope to explore the nonsymmetric and braided versions of this construction
in future work, as well as the connections to topological quantum field theory
\cite{Baez95}.

\begin{defn}{\em
Let M be a commutative monoid with identity e. Define a {\em crossed 
$M$-set} to be a (left) M-set X, together with a function 
$|\,\,\,|:X\rarr M$ such 
that $|mx|=|x|$. (This formula is more complicated in the nonabelian case.
With a nonabelian group, we would require that $|gx|=g^{-1}|x|g$.)

\bigskip

Now define a category {\bf XRel} as follows. 
Objects are crossed M-sets, and maps
are relations $R:X\rarr Y$ such that:

\begin{itemize}
\item $xRy \ \Rightarrow \ mxRmy$
\item $ xRy \Rightarrow |x|=|y|$
\end{itemize}
}\end{defn}

Freyd and Yetter construct 
a category where the objects are functions satisfying
precisely these requirements. They use a nonabelian group and the braiding
is the symmetry adjusted appropriately by the action of G. They then use this
category to develop knot invariants \cite{Freyd89}. In subsequent work, Yetter
uses crossed $G$-sets to construct topological quantum field theories 
\cite{Yetter}. See also \cite{Porter95}.

If X and Y are crossed M-sets, define $X\ox Y$ as cartesian product 
with componentwise action, and $|(x,y)|=|x||y|$. The unit is the one 
element set $I=\{*\}$. Define $|*|=e.$

\begin{thm}
{\bf XRel} is a tensored $*$-category.
\end{thm}

Note that {\bf XRel} is not compact. The counit of the adjunction would
be required to satisfy
$*R(x,x)$ for all $x\in X$, but this would hold if and only if 
$|x|^2=e$. This will be our definition of nuclear object.

Now for all $X,Y$, define ${\cal N}(X,Y)\subseteq Hom(X,Y)$ by:

\begin{center}
$R:X\rarr Y$ is nuclear if and only if  $xRy \Rightarrow |x|^2=|y|^2=e$
\end{center}

\begin{thm}
This defines a nuclear ideal for XRel.
\end{thm}

\section{Distributions as Relations}

In this section, we introduce a generalized category of
relations based on the idea of {\it distributions}.  The guiding intuition
is that composition should be determined by an integral of the
form:

\[ \varphi(x,y);\psi(y,z)=\int\varphi(x,y)\psi(y,z)dy.\] 

The viewpoint here is that the notion of integration generalizes the
existential quantification that appears in the definition of relational
composition.  We will refer to this formula as the ``convolution formula.''
We now introduce a framework in which this makes sense.  A naive approach
is to view $\varphi(x,y)$ and $\psi(y,z)$ as real-valued functions.
However, for such a ``category'' to have identities would require an
equation of the form:

\[\int\varphi(x,y)\delta(y,y')dy=\varphi(x,y')\]

\noindent and similarly for left composition.  The ``function'' playing this 
role is in fact the Dirac $\delta$ which is not a function but a 
{\it generalized function} or {\it distribution} in the sense of Schwartz
\cite{Schwartz57,Treves67,Algwaiz92}.  Unfortunately multiplication of
distributions is not always well-defined.  Formulas like the one above are
sensible only for certain limited kinds of distributions.  In the rest of
this section, we review basic facts about distributions and then develop a
theory of what we call ``tame'' distributions for which the above integral
formula makes sense.

Tame distributions are mentioned in the extant literature (see, for
example, Dieudonn\'e's ``Treatise on Analysis'', volume 7, chapter 23,
sections 9 and 10~\cite{Dieudonne88}), but are not given a name.  

\subsection{Basics of Distributions}

Let \cO\ denote a nonempty open subset of ${\sf R}^{n}$.  Let $\cE{\cO}$
denote the set of $C^{\infty}$ (smooth) functions on \cO\ and $\cD{\cO}$
denote the smooth (complex-valued) functions of compact support on \cO.  
We will refer to the elements of $\cD{\cO}$\ as {\it
test functions}.  
In what follows, we use Greek letters such as $\phi,\psi,\eta$ as test
functions. $\cD{\cO}$\ is given the structure of a 
topological vector space as follows. This structure is described for 
example in \cite{Algwaiz92,Treves67}. 

We begin by considering a compact subset $K\subseteq\cO$, and letting
$\cD{\cO;K}$ be the set of continuous functionals on \cO\ 
with support contained
in $K$. Then we define a family of seminorms on $\cD{\cO;K}$ by the following 
formula,
where $\partial^{x_{1}^{i}\ldots x_{n}^{j}}$ denotes the partial derivative 
with respect to the listed variables:

\[ |\varphi|_{m}=sup\{|\partial^{x_{1}^{i}\ldots x_{n}^{j}}\varphi(x)|\colon
x\in K \mbox{ and } i+\ldots +j\leq m\}\]

We then give $\cD{\cO;K}$ the least topology such that each of these seminorms
is continuous. The existence of such a topology is proved on page 12 of 
\cite{Algwaiz92}. With this topology, $\cD{\cO;K}$ is a {\it Fr\'echet space},
i.e. it is locally convex, metrizable and complete. 

Now observe that

\[\cD{\cO}=\bigcup\{\cD{\cO;K}| ~ K\subseteq\cO \mbox{ and $K$ is compact}\}\]

We then give $\cD{\cO}$ the finest locally convex topology such that the 
inclusions $\cD{\cO;K}\subseteq\cD{\cO}$ are continuous for every compact $K$.
This is known as the {\it inductive limit} of the topologies on $\cD{\cO;K}$.

\begin{thm}\label{continuity}(p.25 \cite{Algwaiz92})
The topology that $\cD{\cO;K}$ inherits as a subspace of $\cD{\cO}$ is the
same as its original topology for every compact $K$. A linear functional
on $\cD{\cO}$ is continuous if and only if the restriction to $\cD{\cO;K}$
is continuous for every compact $K$.
\end{thm}

With this topology, $\cD{\cO}$ is not metrizable. However it is an 
{\it LF space (locally Fr\'echet)} in the sense of \cite{Treves67} p.126.
As such, it is locally convex, Hausdorff and complete.

Then we define a {\it distribution} on
\cO\ to be a continuous, linear (complex-valued) functional on $\cD{\cO}$.
Let ${\cal D}'(\cO)$ denote the set of all distributions on $\cO$. 
Let ${\cal D}'(\cO)$ be given the {weak topology}, p. 45 \cite{Algwaiz92}
or p.197 \cite{Treves67}. This is equivalent to the topology of pointwise
convergence, and ${\cal D}'(\cO)$ is locally convex, Hausdorff and complete.
We will also have need of the following extension theorem \cite{Treves67} p.39.

\begin{thm}\label{extension}
Let $E,F$ be two Hausdorff topological vector spaces, with $A$ a dense
subset of $E$ and $f$ a continuous linear mapping of $A$ into $F$. If $F$ is 
complete, then there is a unique continuous linear mapping $\overline{f}$
from $E$ into $F$ which extends $f$.
\end{thm}

We now describe some examples. 

\begin{enumerate}
\item Let $\cLloc{\cO}$ denote 
the space of locally integrable functions.  Suppose
that $f\in\cLloc{\cO}$.  Define a distribution $T_{f}$ by:

\[ T_{f}(\varphi)=\int_{\cO} f(x)\varphi(x)dx\]

\noindent Note 
that two locally integrable functions determine the same distribution if
and only if they are equal almost everywhere~\cite{Treves67}. A distribution 
of this form is called {\it regular}, and the function $f$ is called the 
{\it kernel of the distribution.}
A distribution
which does not arise in this way 
is called {\it singular}.  Regular distributions are 
fundamental examples, in fact there are a number of strong results
regarding the approximation of distributions by regular
distributions~\cite{Treves67}. 
This justifies thinking of distributions as {\it generalized functions}.

\item As a special case of the previous example, we observe that every
test function is itself locally integrable, and so induces a regular 
distribution.  Thus we have a canonical inclusion 

\[ \iota\colon\cD{X}\hookrightarrow\cDp{X}\]

\noindent given as follows: 

\[ \phi(x)\mapsto [\psi(x)\in\cD{X}\mapsto\int\phi(x)\psi(x)dx]\]

There are similar inclusions for the set of locally 
integrable functions or smooth functions.

\item For any point $x\in \cO$, let $\delta_{x}(\varphi)=\varphi(x)$.
If $0\in \cO$, we denote $\delta_{0}$ simply as $\delta$ and refer to it
as the (one-variable) 
{\it Dirac delta}.  One can show that this distribution is singular,
see for example~\cite{Algwaiz92,Treves67}.  

\item If $\cO\subseteq {\sf R}$,
we may also ``differentiate'' the previous distribution {\it via} the 
formula:

\[ \delta_x'(\varphi)=-\varphi'(x)\]

This distribution is also singular. More generally, if $\cO\subseteq {\sf R}^n$
and $T\in {\cal D}'(\cO)$, we have the formulas:

\[ \frac{\partial}{\partial x_{i}}(T)(\varphi)=-T(\frac{\partial}{\partial 
x_{i}}(\varphi))\]

\[ \partial^{x_{1}^{i}\ldots x_{n}^{j}}(T)(\varphi)=(-1)^{i+\ldots+j}
T(\partial^{x_{1}^{i}\ldots x_{n}^{j}}(\varphi))\]

These formulas allow one to ``differentiate'' nondifferentiable functions, and 
are one of the many advantages of distributions. See, for example, 
\cite{Algwaiz92}, Chapter 2.3.

\item When considering $\cO\times\cO$, we have the {\it trace distribution}, 
\cite{Hormander90} Example 5.2.2, given by:

\[ \varphi\in\cD{\cO\times\cO}\mapsto\int_{\cO}\varphi(x,x)\]

\end{enumerate}

\subsection{The Schwartz kernel theorem}
One is often interested in distributions on product spaces, especially in
the theory of differential equations and their associated Green's functions.  
In this situation the analogy between distributions and
``infinite-dimensional matrices'' is quite striking.  The theory of kernel
distributions can be seen as a formalization of this analogy.  In the
analysis literature,  the notion of ``kernel distribution'' is studied at
length, see for example the massive treatise of
Dieudonn\'e~\cite{Dieudonne88} or the book by Treves~\cite{Treves67}.  
When considering a space of test functions of the form $\cD{X\times Y}$,
there is a canonical subspace of fundamental importance. Consider the tensor
product $\cD{X}\ox\cD{Y}$. A typical element of this space is of the form

\[\sum_{i=1}^{n}\varphi_i\ox\psi_i \mbox{ where $\varphi_i\in\cD{X}$\ and
$\psi_i\in\cD{Y}$}\]

\noindent There is a canonical inclusion of $\cD{X}\ox\cD{Y}$ into  
$\cD{X\times Y}$ given by:

\[\varphi\ox\psi\mapsto[(x,y)\mapsto\varphi(x)\psi(y)]\]

The result we will have use for is:

\begin{prop}
\label{dense}
The space $\cD{X}\ox\cD{Y}$ is sequentially dense in 
$\cD{X\times Y}$.
\end{prop}

Now we have a chance of defining functions on $\cD{X\times Y}$ as the
unique continuous extension of functions defined on $\cD{X}\ox\cD{Y}$ 
using Theorem~\ref{extension}.  

One of the fundamental results in the theory of distributions is the
Schwartz kernel theorem, which gives conditions under which maps from
$\cD{X}$ to $\cDp{Y}$ can be realized as distributions on $X\times Y$.  We
need the following notations to state the theorem.  
If $f$ is a distribution on $X\times Y$ and $\phi\in\cD{X}$ then
$f_*(\phi)$ will be the function from $\cD{Y}$ to the base field given by
$\psi\in\cD{Y}\mapsto f(\phi\ox\psi)$ and
$f^*(\psi)$ is given by the evident ``transpose'' formula.  We have not
yet said that $f_*(\phi)$ and $f^*(\psi)$ are distributions; that is
part of the content of the kernel theorem.  

The Schwartz kernel theorem states:
\begin{thm}
Let $X$ and $Y$ be two open subsets of ${\sf R}^{n}$ and ${\sf R}^{m}$.

\begin{enumerate}
\item Let $f$ be a distribution on $X\times Y$.  For all functions
$\phi\in\cD{X}$ the linear map $f_*(\phi)$ is a distribution on $Y$.  
Furthermore, the map $\phi\mapsto f_*(\phi)$
from $\cD{X}$ to $\cDp{Y}$ is continuous, when  $\cDp{Y}$ is given the
weak topology. 
\item Let $f_*$ be a continuous 
linear map from $\cD{X}$ to $\cDp{Y}$. Then there exists a unique
distribution on $X\times Y$ such that for $\phi\in\cD{X}$ and
$\psi\in\cD{Y}$ the following holds:
\[ f(\phi\ox\psi) = f_*(\phi)(\psi) \]
\end{enumerate}
\end{thm}

Evidently, by symmetry, the same result applies for $f^*$.
In light of the kernel theorem, we may now state the following definition.

\begin{defn}
Suppose that $f$ is a  distribution on $X\times Y$, then we obtain 
the following continuous maps 
(supposing that $\phi\in\cD{X},\psi\in\cD{Y}$ are arbitrary):
\begin{enumerate}
\item $f_*\colon\cD{X}\rarr\cDp{Y}$ is given by 
$f_*(\phi)(\psi)=f(\phi\ox\psi)$
\item $f^*\colon\cD{Y}\rarr\cDp{X}$ is given by
$f^*(\psi)(\phi) =f(\phi\ox\psi)$
\end{enumerate}
\end{defn}

\subsection{Tame Distributions}

To pass from the ``discrete'' category of ordinary relations to a category
of ``continuously varying'' relations, we should replace the usual notion of
morphism in $\rel{}$, a function $X\times Y\rarr 2$, with an integrable
function $X\times Y\rarr \reals{}$, where $X$ and $Y$ are now open subsets
of some Euclidean space.  We have already seen, however, that functions do not
suffice.  One must pass to a class of generalized functions or
distributions.  While distributions satisfy many properties of functions,
they cannot be multiplied and hence the composition formula that we had
proposed does not make sense.  Thus our goal is to introduce a class of
distributions which are sufficiently ``functional'' as to allow us to compose
them using the integral formula discussed above.

We will use a notion defined by Dieudonn\'e in \cite{Dieudonne88}. It will 
provide the first step towards defining a composable class of distributions.
Note that $\cE{X}$ is the space of all smooth complex-valued functions on
$X$ (not necessarily of compact support). Unfortunately, Dieudonn\'e
uses the term {\it regular} which conflicts with the terminology above. 
We therefore use the term {\it Dieudonn\'e-regular}.
  
\begin{defn}{\em
We say that a distribution $f\in\cDp{X\times Y}$ is {\em 
Dieudonn\'e-regular} if
\begin{enumerate}
\item For all functions $\phi\in\cD{X}$, $f_*(\phi)$ is in $\cE{Y}$, that 
is to say there exists $\hat{\phi}\in \cE{Y}$ such that the distribution
$f_*(\phi)\in\cDp{Y}$ is defined by:

\[f_*(\phi)(\psi)=\int_Y\hat{\phi}(y)\psi(y)\]

\item Similarly, for all functions $\psi\in\cD{Y}$, $f^*(\psi)$ is in $\cE{X}$.
\end{enumerate}
}\end{defn}

An equivalent statement is that the function $f_{*}:\cD{X}\rarr\cDp{Y}$
specified by the kernel theorem factors through the inclusion 
$\cE{Y}\hookrightarrow\cDp{Y}$, and similarly for $f^*$.

We would like to define our composition as follows. Given distributions
$f\in \cD{X\times Y}, g\in\cD{Y\times Z}$ which are Dieudonn\'e-regular, 
we try to define a distribution $f;g\in\cD{X\times Z}$ using the following 
formula (with $\phi\in\cD{X},\gamma\in\cD{Z}$).

\[ f;g(\phi\ox\gamma)=\int_Y\hat{\phi}~ ~\hat{\gamma}\]

\noindent Here $\hat{\phi}$ is the element of $\cE{Y}$ associated to the 
distribution $f^*(\phi)$, and $\hat{\gamma}$ is the element of $\cE{Y}$
associated to the distribution $f_*(\gamma)$.

However, the above integral may well be infinite. 
Thus we must add an additional
assumption which assures the finiteness of this integral. One possibility
is to require not only that the two kernels be smooth, but that they have 
compact support.\footnote{In fact, one could use a more 
general class of functions, such as the square integrable 
functions, but we prefer the 
symmetry of the present definition.} Thus, we have the following:

\begin{defn}{\em
A {\em tame distribution} on $X\times Y$ is a distribution $f$ on $X\times Y$
such that each of $f^*$ and $f_*$ factor continuously 
through the appropriate $\iota$, where $\iota$ is the inclusion
of the space of test functions into the space of distributions.
Explicitly, there exist continuous linear maps 

\[f_L\colon\cD{X}\rarr\cD{Y}\]
\[f_R\colon\cD{Y}\rarr\cD{X}\] 

\noindent such that for every $\phi\in\cD{X}$ and \ $\psi\in\cD{Y}$, we have:

\[f_*(\phi)(\psi) = f^*(\psi)(\phi) =
f(\phi\ox\psi) =\int f_L(\phi)\psi dy = \int\phi f_R(\psi) dx\]

}\end{defn}

Note that we are not saying that $f_L$ and $f_R$ have functional kernels
and certainly not that $f$ has a functional kernel.  But rather that $f^*$
and its adjoint $f_*$ map test functions to distributions with test
functions as kernels.  In some sense, tame distributions are allowed to be
mildly singular, in that composing with a
test function ``tames'' the singularity.

Dieudonn\'e, in \cite{Dieudonne88}, page 77, examines the question of when
the operators $f^*$ and $f_*$ map test functions to test functions, and he 
derives the following theorem.

\begin{thm}
Let $f$ be a Dieudonn\'e-regular distribution on $X\times Y$.
The following are equivalent:
\begin{enumerate}
\item The operator $f_*$ extends to a continuous linear map from the
Fr\'echet space $\cE{X}$ to the Fr\'echet space $\cE{Y}$.
\item The operator $f^*$ maps $\cD{Y}$ to $\cD{X}$.
\item The operator $f^*$ maps $\cEp{Y}$ to $\cEp{X}$, where $\cEp{Y}$ is the
space of distributions of compact support (see \cite{Algwaiz92} for 
the definition of support of a distribution).
\end{enumerate}
\end{thm}

\subsection{Examples}

\begin{itemize}
\item Let $X$
be an open subset of $\reals^n$.  The trace distribution on $X\times X$
is given by $Tr(\eta) =\int \eta(x,x) dx$ where
$\eta(x,x')\in\cD{X\times X}$.  From this definition it follows that
$Tr_*(\phi)(\psi) = Tr^*(\psi)(\phi) = Tr(\phi\ox\psi) =
\int \phi(x)\psi(x) dx.$  Thus we clearly have $Tr_L(\phi)
=Tr_R(\phi) =\phi$, which shows that $\delta$ is tame.  This tame 
distribution will act as the identity in our category.

\item Suppose that $T$ is a regular distribution on $X\times Y$ with a test
function $\beta(x,y)$ as its kernel, that is to say:

\[ T(\alpha(x,y))=\int_{X\times Y}\beta(x,y)\alpha(x,y)\]

\noindent Then $T$ is tame with its associated functions being given by:

\[ T_{L}(\phi)=\int_X\beta(x,y) \phi(x)\]
\[ T_{R}(\psi)=\int_Y\beta(x,y) \psi(y)\]

We write $\cT{X,Y}{}$ for the tame distributions on $X\times Y$.  
 
\end{itemize}

\subsection{Composing tame distributions}

Given tame distributions we can define the following operation which will
serve as composition.  Suppose that $f\in\cT{X,Y},g\in\cT{Y,Z}$. 
We define $f;g\in\cT{X,Z}$ as follows.
Given that $f$ is tame, we have a continuous function $f_L:\cD{X}\rarr\cD{Y}$.
Applying the first part of the 
Schwartz kernel theorem to $g$, we obtain a morphism
$g_*:\cD{Y}\rarr \cDp{Z}$. Composition gives a continuous map
$\cD{X}\rarr\cDp{Z}$. By the second part of the kernel theorem, we obtain
a distribution on $X\times Z$.

Alternatively, we could use the extension theorem, Theorem \ref{extension}.
Let $\phi\in\cD{X},\psi\in\cD{Z}$.  We set
\[ (f;g)(\phi\ox\psi) = \int f_L(\phi)g_R(\psi) dy.\]
This, of course, only defines $f;g$ on $\cD{X}\ox\cD{Z}$ rather than on
$\cD{X\times Z}$. We then use the fact that
the tensor product is a dense subspace to extend composition to all of 
$\cD{X\times Z}$.
One observes that $f;g$ is tame
as can be seen by an elementary calculation, noting $(f;g)_L = f_L;g_L$
and $(f;g)_R = g_R;f_R$ and the tameness of $f$ and $g$.

\subsection{The category \drel}

\begin{defn}{\em
The category \drel\ has as objects open subsets on $\reals^n$, 
and, as morphisms, tame distributions.  Composition is as described above.}
\end{defn}

\begin{thm}
\drel\ is a tensored $*$-category.
\end{thm}
\begin{proof}{}
Evidently we can verify properties of
the composition $f;g$ by carrying out calculations on the distribution defined
on $\cD{X}\ox\cD{Z}$ and appealing to continuity and the density of
$\cD{X}\ox\cD{Z}$ in $\cD{X\times Z}$.  We have already noted above that
$f;g$ is tame.  A simple calculation shows that the trace distribution is
the identity for composition.  

To verify associativity we calculate as follows.  Let $f\in\cD{X\times
Y}$, $g\in\cD{Y\times Z}$ and $h\in\cD{Z\times W}$ be tame
distributions.  Then we have:

\begin{eqnarray*}
((f;g);h)(\phi(x)\ox\rho(w))&=& \int (f;g)_L(\phi)h_R(\rho)dz\\
&=& \int g_L(f_L(\phi))h_R(\rho)dz\\
&=& \int f_L(\phi) g_R(h_R(\rho)) dy\\
&=&(f;(g;h))(\phi\ox\rho)
\end{eqnarray*}

Thus we have shown that \drel{} is a category.  The tensor product is
given as follows.  Given objects $X$ and $Y$ we define $X\ox Y$ as the
cartesian product space $X\times Y$.  Given morphisms in \drel{}
$f\colon X\rarr Y$ and $g\colon X'\rarr Y'$ we can 
define $f\ox g\colon X\ox X'\rarr Y\ox
Y'$ as follows.  We first define $f\ox g$ as a distribution on
$\cD{X}\ox\cD{X'}\ox\cD{Y}\ox\cD{Y'}$ by the formula
\mbox{$(f\ox g)(\phi(x)\ox\phi'(x')\ox\psi(y)\ox\psi'(y')) =
f(\phi\ox\psi)g(\phi'\ox\psi')$}.  It is routine to verify that this is
tame.  We extend $f\ox g$ to all of $\cD{X\times X'\times Y\times Y'}$
as above.  The one-point space, written $I=\{*\}$, is the unit
for the tensor (with measure $\mu(\{*\})=1$).

Finally the $*$-structure is the identity on objects.  On morphisms, the
only thing that changes is the role of $f_L$ and $f_R$. The conjugate functor 
is taken to be the identity.

\end{proof}

\begin{rem}
\label{homi}
As an example, we will describe $Hom(I,X)$, where $X$ is an arbitrary object. 
Clearly, $\cD{I}$ is isomorphic to the base field. We must have two functions:

\[f_L\colon\cD{I}\rarr \cD{X}\]
\[f_R\colon\cD{X}\rarr \cD{I}\]

\noindent such that, for all $\psi\in\cD{X}$:

\[\int_X f_L(1)\psi=\int_I 1 f_R(\psi)\]

But evidently $\int_I 1 f_R(\psi)=f_R(\psi)$. So the function $f_R$ is uniquely
determined by the function $f_L$. Hence we may conclude that 
$Hom(I,X)$ is in bijective correspondence to test functions on $X$.
\end{rem}

We now display a nuclear ideal for \drel{}.  We remarked that
not all tame distributions can be viewed as integral operators with
functions as kernels.  In particular the identity morphisms do not have
this property.  However, we will see that tame distributions with 
functional kernels form a nuclear ideal.
\begin{defn}{\em
Given objects $Y$ and $Z$ of \drel{} we define the set of {\em nuclear
morphisms}, written \nuc{Y}{Z}, as the collection of tame distributions
$g\colon Y\rarr Z$ such that 
$\exists\beta(y,z)\in\cD{Y\times Z}$ 
with the property that for every $\phi(y,z)\in\cD{Y\times Z}$:

\[g(\phi)=\int \beta(y,z)\phi(y,z)dydz\]  
}\end{defn}

Note that the test function $\beta(y,z)\in\cD{Y\times Z}$ associated to
to the tame distribution $g$ is unique. Thus, the set \nuc{Y}{Z} is in 
bijective correspondence to $\cD{Y\times Z}$.

\begin{thm}\label{drelnuc}
The sets \nuc{Y}{Z} form a nuclear ideal for \drel{}.
\end{thm}
\begin{proof}{}
As already remarked, if $g\in\nuc{Y}{Z}$
and if $\beta$ is its kernel, then: 

\[\forall\psi\in\cD{Y}, ~ g_L(\psi)=\int \beta(y,z)\psi(y)dy\]

To verify that we have an ideal, we have to show that for
any $f\in\cT{X,Y}$ the composite $f;g$ is nuclear and symmetrically
for composition on the other side of $g$.  In order to verify this we
need to find a kernel for $f;g$.  We claim that this kernel is
$\alpha(x,z) =_{df} f_R(\beta(y,z))$ where we interpret this formula as
follows.  For each fixed $z\in Z$ $\beta(y,z)$ is a smooth function of compact
support in $Y$; $f_R$ acts on this function to produce a function of
compact support in $X$. The function $\alpha(x,z)$ evidently has compact
support, and its smoothness is a consequence of the continuity 
of $f_R$. It suffices to prove this for functions $\beta$ of 
the form $\beta(y,z)=a(y)b(z)$ where
$a\in\cD{Y}$ and $b\in\cD{Z}$. This follows from Proposition~\ref{dense}
which implies that arbitrary $\beta$ can be written:

\[ \beta(y,z)=\lim_{n\rarr\infty}\sum_{i=1}^{m_n}a_{i,n}(y)b_{i,n}(z)\]

The general result then follows from the linearity and continuity of $f_R$.

\noindent Now observe that for a fixed $z$:

\[
f_R(\beta(y,z))=f_R(a(y)b(z))
=f_R(a)(x)b(z)
\]

Now we calculate as follows, again letting $\beta(y,z)=a(y)b(z)$ and relying
on linearity and continuity for the general result:

\begin{eqnarray*}
(f;g)(\phi(x)\ox\psi(z)) &=& \int_Y f_L(\phi)(y)g_R(\psi)(y) dy\\
&=& \int_Y f_L(\phi)(y)[\int_Z \beta(y,z)\psi(z)dz] dy\\
&=&\int_Z [\int_Y f_L(\phi)(y)\beta(y,z) dy] \psi(z) dz\\
&=&\int_Z [\int_Y f_L(\phi)(y) a(y)b(z) dy]
\psi(z) dz\\ 
&=&\int_Z [\int_Y f_L(\phi)(y)a(y) dy]
b(z)\psi(z)dz\\ 
&=&\int_Z[\int_X\phi(x)f_R(a)(x)dx]b(z)\psi(z)dz\\
&=&\int_X\phi(x)[\int_Z
[f_R(a)(x)b(z)]\psi(z)dz]dx\\ 
&=& \int_X\int_Z\phi(x)\alpha(x,z)\psi(z) dxdz.
\end{eqnarray*}

It follows that $f;g$ is an
integral operator with $\alpha$ as its kernel.  The verification for
composition on the other side of $g$ is very similar.

To complete the proof, we need to show that $Hom(I,X\ox
Y)\cong\nuc{X}{Y}$.  This isomorphism is described in Remark~\ref{homi}.
It remains to verify the equations. Naturality requires an argument 
similar to the previous calculation. Compactness is quite straightforward.

\end{proof}

\section{The Category {\bf PRel}}

In this section, we define a category of \emph{probabilistic relations},
and describe a nuclear ideal for it.  We will see that we indeed get most
of the important properties of the category of relations, i.e.  we have a
tensored $*$-category with a nuclear ideal.  Thus one may think of this
category as representing relations ``smeared out probabilistically''.  Once
again, as in \drel\, we have a situation where the identity maps are too
singular to be in the nuclear ideal.  The nuclear ideal can be thought of
as functions but the ambient category has to be described in terms of measures.

\subsection{Basic Definitions of Measure Theory} We assume the reader is
familiar with the basic concepts of measure theory.  We recall the basic
definitions for completeness.  A reader who remembers these definitions can
skip to the start of the next section without loss of continuity.

\begin{defn} {\em A $\sigma${\em -field} $\Sigma$ on a set $X$ is a
collection of subsets of $X$ which 

\begin{enumerate} 
\item includes the whole space $X$, 
\item is closed under complementation, and 
\item is closed under finite and countable unions.  
\end{enumerate} 
A {\em measurable space} is a set 
together with a $\sigma$-field.  A {\em measurable
function} from a measurable space $(X,\Sigma_X)$ to $(Y,\Sigma_Y)$ is a
function from $X$ to $Y$ such that for all $B\in\Sigma_Y$ we have
$f^{-1}(B)\in\Sigma_X$.}  
\end{defn} 

Given a measurable space
$(X,\Sigma_X)$, 
we call the members of $\Sigma_X$ {\em measurable sets}.  If $B$
is a measurable set then the characteristic function of $B$ is
denoted $\chi_B$ and is clearly measurable.  

\begin{defn} {\em A
{\em measure} $\mu$ on a measurable space $(X,\Sigma_X)$ is a function
$\mu:\Sigma_X\to [0,\infty]$ such that 

\begin{enumerate} 
\item $\mu(\emptyset) = 0$ 
\item if $\{A_i|i\in I\}$ is a pairwise-disjoint
family of measurable sets, with $I$ countable, then 

\[ \mu(\cup_{i\in I}A_i)= \sum_{i\in I}\mu(A_i). \] 
\end{enumerate} 

If we have a measure taking
values in $[0,1]$ we call it a {\em sub-probability measure} and if the
measure (``mass'') of the whole space is $1$ we say that it is a
{\em probability measure}.  A $\sigma$-field equipped with a measure is
called a {\em measure space} and equipped with a probability measure it
is called a {\em probability space}.  
}\end{defn}

Sets of measure zero play an important role.  The phrase {\em almost
everywhere} is frequently used to assert that a certain property holds
everywhere except on a set of measure zero.  If there is confusion about
which measure is intended we might say, for example, $P$-almost everywhere.

The set of real numbers and the closed unit interval
$[0,1]$ play a central role in the subsequent discussion.  As
measurable spaces, each has two $\sigma$-fields which are often used, the
\emph{Borel} $\sigma$-field and the \emph{Lebesgue} $\sigma$-field.  Any
collection of subsets of a set $X$ generates a $\sigma$-field,
namely the least $\sigma$-field containing all the sets of the given
collection.  If we take the open sets of any topological space and generate
a $\sigma$-field we get the Borel $\sigma$-field.  In particular we get the
Borel $\sigma$-field on the reals.  This $\sigma$-field on the reals can be
given a measure in such a way that the measure of an interval is its
length.  The resulting measure space has the property that there are
subsets of sets of measure zero that are not measurable.  There is a
canonical ``completion'' procedure which yields an extended $\sigma$-field
and measure, such that any previously measurable set has the same measure
and all subsets of sets of measure $0$ are measurable (and have measure
$0$).  When applied to the Borel subsets of the reals with the Lebesgue
measure one gets Lebesgue measurable sets (with the Lebesgue measure).
\emph{In our discussion we always mean Borel measurable whenever we talk
about a measurable subset of the reals.}  

In some older books~\cite{Halmos50,Rudin66}, a measurable function from the
reals to the reals is defined to be a function where the inverse image of
an open set has to be a Lebesgue measurable set rather than a Borel
measurable set.  This has the unfortunate effect that the composite of two
measurable functions need not be measurable.
A suitable reference for the above discussion is \cite{Malliavin95},
but any good book on probability theory such as Ash~\cite{Ash72},
Billingsley~\cite{Billingsley95}, or Dudley~\cite{Dudley89} covers this 
material.

\subsection{A category of stochastic kernels}

Probability theory has been examined in the past from a categorical
perspective.  For example, Giry~\cite{Giry80} has given the following
construction, based on hints in unpublished notes of Lawvere.  Wendt has
examined this construction extensively \cite{Wendt93,Wendt94}.

Let \meas\ denote the category of measurable spaces and measurable
functions.  We will now describe a triple $T$ on the category \meas.  In
what follows, when we talk about measurable functions into $[0,1]$, we
always mean the Borel $\sigma$-field on $[0,1]$, denoted ${\cal B}$.  If
$(X,\Sigma)$ is an object of \meas, then we define $T(X,\Sigma)$\ to be the
set of probability measures on $(X,\Sigma)$\ equipped with the least
$\sigma$-algebra making the evaluations

\[ e_{B}\colon T(X)\rarr [0,1] \mbox{ defined by } e_{B}(P)=P(B)\]

\noindent measurable, where $B$\ ranges over the measurable sets of $X$.
$T$\ acts on maps by the formula:

\[ T(f)(P)(B')=P(f^{-1}(B'))\]

\noindent where $f\colon X\rarr Y$ \ and $B'\in \Sigma_{Y}$.

The unit for the triple $\eta\colon id\rarr T$\ is defined by the formula:

\[ \eta_{X}(x)(B)=\chi_{B}(x)\]

\noindent where $x\in X$ \ and $\chi_{B}$ is the characteristic function of
$B$.

The multiplication $\mu\colon T^{2}\rarr T$
is defined as follows.  If $P'\in T^{2}(X)$, then $P'$\ defines a measure on
$T(X)$, and we use it to form the following integral:

\[ \mu_{X}(P')(B)=\int_{T(X)} e_{B}dP'\]

With these definitions, one can then prove~\cite{Giry80}:

\begin{thm}
$(T,\eta,\mu)$ \ form a triple on \meas.
\end{thm}
 
To understand the structure of the Kleisli category, we require the following 
definition.

\begin{defn}{\em 
If $(X,\Sigma)$\ and $(X',\Sigma')$ \ are measurable spaces, then 
a} stochastic kernel  {\em on $X\times X'$ \ is a function

\[ \rho\colon X\times \Sigma'\rarr [0,1]\]

\noindent that is measurable in its first argument, for each fixed
measurable set and 
a probability measure in its second argument for each point in $X$.}
\end{defn} 
Stochastic kernels are closely related to regular conditional probability
distributions~\cite{Ash72,Dudley89}.

\bigskip

If $\tau$ is a stochastic kernel on $X\times Y$ 
and $\rho$\ is a stochastic kernel on $Y\times Z$, 
then we can compose $\rho$ \ and $\tau$ to obtain a 
stochastic kernel $\tau\circ\rho\colon X\times \Sigma_{Z}\rarr [0,1]$,
using the following formula:

\[\tau\circ\rho(x,C)\colon 
\int_{Y}\rho(-,C)d\tau(x,-) \mbox{ for all $x\in X,
C\in \Sigma_{Z}$}\]

\noindent Note that in the above formula $\rho(-,C)$\ is acting as the 
measurable function, and $\tau(x,-)$\ as the measure.  The associativity
of this composition follows easily from the monotone convergence theorem.

So we obtain a category \cpd, whose objects are measurable spaces, and whose
morphisms are stochastic kernels.  The identity for this
category is given by the $\delta$-formula:

\[ \delta(x,A) = \left\{ \begin{array}{ll}
                          1 & \mbox{if $x\in A$} \\
                          0 & \mbox{otherwise}
                         \end{array}
                 \right.  \]

\bigskip

One can now derive \cite{Giry80}:

\begin{thm} 
The Kleisli category for the triple $T$ \ is equivalent to \cpd.
\end{thm}

Given a morphism $f\colon X\rarr TY$ in the Kleisli category, one obtains
a stochastic kernel \emph{via} the formula:

\[ F\colon X\times\Sigma_{Y}\rarr[0,1] \mbox{ is defined by } F(x,B')=f(x)(B')\]

\subsection{Probabilistic Relations}
While the category \cpd\ allows valuable insights into probability theory -
for example, the Chapman-Kolmogorov equation is simply functoriality
\cite{Giry80} - it lacks some of the structure one requires of a category
of relations; notably the ability to take the converse.  To pass to a
category which is more relational in nature, we will use measures on the
product space.  Unfortunately one cannot compose measures in any simple
way.  Given measures on the product space, there is no obvious sense in
which one can integrate them to compose as in the category \cpd{}.  The
idea is to rely on a basic theorem which says that given such product
measures, \emph{on suitable spaces}, one can construct a pair of stochastic
kernels -- which, together with the marginal distributions, determine the
original measure on the product space -- and then compose them in the
manner described for \cpd{}.

We now give the details of the construction.  First suppose that we have a
pair of measurable spaces $(X,\Sigma_X)$ and $(Y,\Sigma_Y)$, a probability
measure $P_X$ on $(X,\Sigma_X)$, and a stochastic kernel $h(x,B)\colon
X\times \Sigma_Y\rarr [0,1]$.  Then we have a \emph{unique} measure $P$ on
the product such that for all $A\in \Sigma_X$:
\[ P(A\x B) = \int_A h(x,B)dP_X(x).\]
Thus if we have a pair of stochastic kernels $h:X\x\Sigma_Y\to [0,1]$ and
$k:Y\x\Sigma_X\to [0,1]$ and probability distributions $P_X$
on $(X,\Sigma_X)$ and $P_Y$ on $(Y,\Sigma_Y)$ -- satisfying an evident
compatibility condition -- we can reconstruct a unique probability measure
on the product space.  

Conversely, given a measure $P$ on the product $X\x Y$ we can construct a
measure on each of the factor spaces by setting $P_X(A) := P(A\x Y)$ and
$P_Y(B):= P(X\x B)$.  These are called the \emph{marginals}.  Knowing one
of the marginals and the appropriate stochastic kernel is equivalent to
knowing the product measure.  Clearly the pair of stochastic kernels does
not uniquely determine the product measure; it does not even determine the
marginals.  We now need to show how to go from the product measure to the
stochastic kernels.

The situation we have is: a pair of measure spaces $(X,\Sigma_X,\mu_X)$ and
$(Y,\Sigma_Y,\mu_Y)$ and a measure, say $\alpha$, on the product space
equipped with the product $\sigma$-field, $\Sigma_X\ox\Sigma_Y$.  We want
to construct a stochastic kernel $h:X\x\Sigma_Y\to[0,1]$.  The product
space is a product in the category \meas\, and is equipped with the usual
projections $\pi_1$ and $\pi_2$ to $X$ and $Y$ respectively.  We want to
construct $h:X\x\Sigma_Y\to [0,1]$ as in the diagram
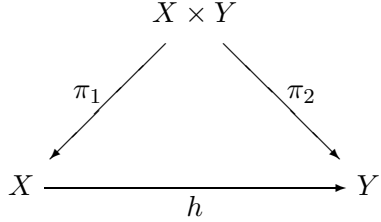
\[
\begin{diagram}
 & & X\x Y & & \\
 &\ldTo^{\pi_1} & & \rdTo^{\pi_2}& \\
X & &\rTo_{h} & & Y
\end{diagram}
\]
such that 
\[ \int_A h(x,B)\mu_X =  \alpha(A\x B). \]
where $h$ is the morphism (of the category \cpd{}) that we are trying to
construct and $\pi_1,\pi_2$, the projections, are morphisms of the category
\meas.  However, this construction requires some assumption on the spaces
involved.  

More precisely, we require that the spaces are \emph{Polish
spaces}\footnote{We could have more general spaces, for example analytic
spaces \cite{Hoffman-Jorgenson94}.}.  Recall that a Polish space is the topological space underlying a
complete separable metric space.  This assumption is quite common in
probability theory and allows the construction of \emph{regular conditional
probability distributions}~\cite{Ash72,Billingsley95,Dudley89}.  We will not
invoke these general concepts here.

We state a slightly more general theorem from which the construction of $h$
in the preceding paragraph follows immediately.  
\begin{thm}\label{bet}
Suppose that $(U,\Sigma_U,P)$ is a probability space, $V$ is a Polish space
with the Borel $\sigma$-field, written $\Sigma_V$, and $(W,\Sigma_W)$ is a
measurable space.  Suppose that $f$ is a measurable function from $U$ to
$V$ and that $g$ is a measurable function from $U$ to $W$. Then there exists
a \cpd{} morphism, i.e. a stochastic kernel, $Q:W\to V$ as shown in the
diagram 
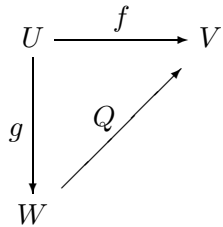
\[
\begin{diagram}
U & \rTo^f & V \\
\dTo^g & \ruTo^Q & \\
W & & 
\end{diagram}
\]
such that for all $A\in\Sigma_W,B\in\Sigma_V$:
\[ \int_{g^{-1}(A)} Q(g(u),B)\mbox{d}P(u) 
= P(g^{-1}(A)\cap f^{-1}(B)).\]
This $Q$ is unique in the sense that if $Q'$ is another stochastic kernel
satisfying the same equation then for $P$-almost all $u\in U$ $Q(u,\cdot)$
and $Q'(U,\cdot)$ are identical.
\end{thm}
Roughly speaking, this says that $Q$ composed with $g$ agrees with $f$ at
least when evaluated on the measures $P$.  
In probability texts this theorem is stated in terms of existence of
regular conditional probability distributions \emph{relative to a sub
$\sigma$-field}.  We have essentially the same situation since the set of
inverse
images under $g$ of the $W$-measurable sets forms a sub-$\sigma$-field of
$\Sigma_U$.  With this identification, theorem~\ref{bet} is equivalent to
theorem 10.2.2 of ~\cite{Dudley89}.  

We are now ready for the corollary of chief interest.

\begin{cor}\label{prod}
Given Polish spaces $X$ and $Y$ with their Borel $\sigma$-fields and a
probability measure $\alpha$ on the product space, there is a stochastic
kernel $Q_1(x,B)$ (i.e. a \cpd{} morphism from $X$ to $Y$), where
$B\in\Sigma_Y$ and a stochastic kernel  $Q_2(y,A)$ (i.e. a \cpd{} morphism
from $Y$ to $X$), where $A\in\Sigma_X$, such that 
\[\int_{A}Q_1(x,B)d\alpha_{X}=\alpha(A\times B) =
\int_{B}Q_2(y,A)d\alpha_{Y}. \]
\end{cor}
\begin{proof}{}
We use the theorem~\ref{bet} with $X\x Y$ as $U$, $X$ as $W$ and $Y$ as $V$
and the projection maps as $f$ and $g$.  Now we immediately get $Q_1$.
To see that the equation is satisfied we check as follows:
\[ \alpha(\pi_1^{-1}(A)\cap\pi_2^{-1}(B))  = \alpha(A\x B).\]
On the other hand the left hand side of the equation asserted in
theorem~\ref{bet} is, in this case,
\[\int_{A\x Y}Q_1(\pi_1(\langle x,y \rangle),B)d\alpha.\]
This can be rewritten as
\[\int_A Q_1(x,B)d\alpha\circ\pi_1^{-1} = \int_A
Q_1(x,B)d\alpha_X\]
which is the desired result.  One gets the result for $Q_2$ similarly.
\end{proof}

Here are two simple example applications of corollary~\ref{prod}.  For the
first we take the product measure $\alpha$ to be $\mu\ox\nu$.  In this case
the stochastic kernel $h:X\x\Sigma_Y\rightarrow [0,1]$ is $h(x,B) =
\nu(B)$, i.e. it is independent of $x$.  If we take the product $X\x X$
with the measure $\Delta$ defined by $\Delta(A\x B)=\mu(A\cap B)$, we get
the usual Dirac delta $\delta(x,A)$.  

Finally, to define morphisms in our category, we proceed as follows.
Given two measures, $\mu$ and $\nu$, on a measurable space we say $\nu$
is \emph{absolutely continuous} with respect $\mu$, written $\nu << \mu$,
if for any measurable set $A$, $\mu(A)=0$ implies that $\nu(A)=0$.  We now
assume that the marginal $\alpha_{X}$ is absolutely continuous with respect
to $\mu$.  By applying the Radon-Nikodym theorem~\cite{Billingsley95}, we
obtain a measurable function $h(x)\colon X\rarr {\sf R}$\ such that 
\[\alpha_{X}(A)=\int_{A}h(x)d\mu(x)\]

\noindent From which it follows that:

\[ \int_{A}Q(x,B)d\alpha_{X}(x)=\int_{A}Q(x,B)h(x)d\mu(x)\]

\bigskip

We refer to the function $F(x,B)=Q(x,B)h(x)$\ as the \emph{stochastic
kernel associated to $\alpha$.}   

\begin{defn} {\em We define a category \prel{} as follows.  The objects of 
\prel{} are triples $(X,\Sigma,\mu)$, where $X$ is a Polish space, 
$\Sigma$ the associated $\sigma$ field and 
$\mu$ is a probability measure on $(X,\Sigma)$.  A morphism 
$\alpha\colon(X,\Sigma,\mu)\rarr(X',\Sigma',\mu')$
is a probability measure on $\Sigma\ox \Sigma'$ whose marginals are
absolutely continuous with respect to $\mu$\ and $\mu'$.  }\end{defn}

To compose morphisms $\alpha\colon X\rarr Y$\ and $\beta\colon Y\rarr Z$,
we calculate their associated stochastic kernels  $F(x,B)$\ and
$G(y,C)$ and compose as in the above Kleisli category to obtain a
stochastic kernel $H(x,C)$.  We then obtain a measure on $X\times Z$
\emph{via} the formula:

\[ \gamma(A\times C)=\int_{A}H(x,C)d\mu(x)\]

\begin{thm} \prel{} is a category.
\end{thm}

\begin{prf}{}
The only thing remaining to consider is the identity.  If $(X,\Sigma,\mu)$
is an object, its identity is given by $\Delta(A\times A')=\mu(A\cap A')$,
with the associated conditional distribution given by the Dirac
$\delta$.
\end{prf}

\begin{thm}
\prel{} is a tensored $*$-category.
\end{thm}

\begin{prf}{}
The $*$-structure of \prel{} is evident, and the tensor product 
on objects is given by the product in the category \meas, that is, 
one takes the product of the 2 sets, the tensor of the $\sigma$-algebras,
and the product measure.  The necessary equations are all straightforward 
to verify.
\end{prf}

It is worth understanding the nature of isomorphisms in \prel{} in order to
get a better sense of the role of the measures on the \prel{} objects.  We
consider first objects with the same underlying Polish space and hence
$\sigma$-field.  We will show that two such objects are isomorphic exactly
when they define the same ideal of sets of measure zero.  
\begin{prop}\label{iso1}
Consider two \prel{} objects $X_1$ and $X_2$ where $X_1 = (X,\Sigma,\mu)$
and $X_2 = (X,\Sigma,\nu)$.  They are isomorphic in \prel{} if and only if 
$\mu <<
\nu$ and $\nu << \mu$.
\end{prop}
\begin{prf}{}
Suppose first that $\mu << \nu$ and $\nu << \mu$.  We define an isomorphism
$H:X_1\to X_2$ and $K:X_2\to X_1$ as follows\footnote{As usual we define
measures on product spaces by specifying them on the semi-ring of
``rectangles'' and then relying on the standard extension
theorems~\cite{Billingsley95} to obtain the unique extension to the whole
space.}.  We set $H(A\x B) = \mu(A\cap B)$ and $K(A\x B) = \nu(A\cap B)$.
The marginals are
\[ H_1 = H_2 = \mu \mbox{ and } K_1 = K_2 = \nu.\]
By the absolute continuity assumptions these are \prel{} morphisms.  The
associated stochastic kernels are just the Dirac delta distributions and the
composite of these distributions are again Dirac delta distributions.
As we have observed before the Dirac delta distribution is the stochastic
kernel associated with the identity morphism.  Thus $H$ and $K$ form an
isomorphism.

Conversely, suppose that we have an isomorphism $H:X_1\to X_2$ and
$K:X_2\to X_1$.  Suppose that $\mu(A)=0$ for some $A\in\Sigma$.  Let
$h'$ be the stochastic kernel from $X_2$ to $X_1$ associated with $H$, then
we have 
\[ \int_{X_2} h'(x,A)d\nu(x) = H(A\x X) = H_1(A) = 0\]
where the last equality follows from $H_1 << \mu$ as required for $H$ to be
a \prel{} morphism.  We are writing integrals over $X_2$ and $X_1$ rather
than over $X$ in order to avoid confusion; of course $X_1$ and $X_2$ are
both $X$ as sets.  Since $h'$ is always nonnegative we have that it is
$\nu$-almost everywhere $0$.  Let $k$ be the stochastic kernel from $X_1$
to $X_2$ associated to $K$.  Since $H$ and $K$ form an isomorphism, we have
\[ \int_{X_2} h'(x',A)k(x,dx') = \delta(x,A). \]

Integrating both sides of this equation over $X_1$ using $\nu$, we get
\[ \int_{X_1} [\int_{X_2} h'(x',A)k(x,dx')] d\nu(x) =
\int_{X_1}\delta(x,A) d\nu(x) = \nu(A)\]
It can easily be shown, using the monotone convergence theorem, that we can
rewrite the left hand side as
\[ \int_{X_2} h'(x',A) [\int_{X_1} k(x,dx')d\nu(x)] \]
where the integral in square brackets defines the measure used for the
outer integration.  This measure is absolutely continuous with respect to
$\nu$ since it is defined by $k$.  Since the integrand $h'(x',A)$ is
$\nu$-almost everywhere $0$, the whole integral is $0$.  Thus $\nu(A)=0$ and
$\nu << \mu$.  Similarly $\mu << \nu$.
\end{prf}
\begin{obs}\label{iso2}
{\rm 
Similarly, given two Polish spaces and a Borel isomorphism between them,
one can show that the two objects are isomorphic if and only if the Borel
isomorphism preserves and reflects sets of measure zero.}
\end{obs}

In view of proposition~\ref{iso1} and observation~\ref{iso2} the following
important theorem of classical measure theory (see, for
example, theorem 13.1.1 in~\cite{Dudley89}) almost completes the analysis
of isomorphisms in \prel{}.
\begin{thm}\label{bigth}
If $X$ and $Y$ are Polish spaces, then $X$ and $Y$ are Borel isomorphic if
and only if $X$ and $Y$ have the same cardinality.  Moreover this
cardinality must be either finite, countable or the cardinality of the
continuum. 
\end{thm}
Now we can state the main theorem about isomorphisms in \prel{}.
\begin{thm}
Let $X$ and $Y$ be two objects in \prel{}.  Then $X$ and $Y$ are isomorphic
if and only if there is a Borel isomorphism between them and that
isomorphism preserves and reflects sets of measure $0$.
\end{thm}
\begin{prf}{}
In view of theorem~\ref{bigth}, 
it remains to show that isomorphic objects in \prel{} always have
the same cardinality.  First note that for finite or countable objects in
\prel{} the stochastic kernels are just stochastic matrices.  Thus an
elementary rank argument suffices.  

In the case that one of the objects has an uncountable underlying set we
argue as follows.  It is easy to see that in an uncountable set, with any
$\sigma$-field and with any probability measure, say $P$, there can be at
most countably many points, $x$, with $P(\{x\})\neq 0$.

Now suppose that $(X,\Sigma_X,\mu)$, with $X$ a countable set, and
$(Y,\Sigma_Y,\nu)$, with $Y$ uncountable, are \prel{} objects.  Suppose,
for the moment, that $\mu(\{x\})$ is nonzero for every $x\in X$.  Now
suppose that we have an isomorphism $H:X\to Y$ with inverse $K:Y\to X$.
Thus we have stochastic kernels as follows:
\[ h^+,k^-:X\x\Sigma_Y\to [0,1] \mbox{ and } h^-,k^+:Y\x\Sigma_X\to [0,1].\]
Since these are isomorphisms, we obtain the equation
\[ \int_X k^-(x,B)h^-(y,dx) = \delta(y,B).\]
Since $X$ is countable, this reduces to
\[ \sum_{x\in X} k^-(x,B)h^-(y,\{x\}) = \delta(y,B).\]
Let $B=\{y\}$, where $\{y\}$ is a set with $\nu$-measure zero.  Now observe
that $k^-$ must satisfy
\[\int_X k^-(x,\{y\})d\mu = K(X\x \{y\}) = K_Y(\{y\}) = 0\]
where the last equality is a consequence of the absolute continuity
requirement.  But
\[\int_X k^-(x,\{y\})d\mu=\sum_{x\in X}k^-(x,\{y\})\mu(\{x\}).\]
By assumption, for every $x\in X$, we have that $\mu(\{x\})\neq 0$. 
Thus, for every $x\in X$, it is the case that $k^-(x,\{y\}) = 0$.
So we conclude:
\[ \sum_{x\in X} k^-(x,\{y\})h^-(y,\{x\}) = 0 \neq \delta(y,\{y\}).\]
This is a contradiction.  

Finally, recall that the stochastic kernels are uniquely defined only {\em
almost} everywhere.  In particular, for a countable probability space,
the set of all points of measure zero itself has measure zero. Thus, 
at points where $\mu(\{x\})=0$, we can
define $k^-(x,B)$ to be $0$, and the above argument still applies.
\end{prf}

\subsection{A nuclear ideal for \prel{}}

To determine a nuclear ideal for \prel{}, we must consider the set
$Hom(I,X\ox Y)$.  By definition, this consists of measures $\alpha$\ which
are absolutely continuous with respect to the product measure
$\mu\times\mu'$.  By Radon-Nikodym, we can construct a measurable function
$f\colon X\times Y\rarr {\sf R}$\ such that for all $C\in
\Sigma_{X}\ox\Sigma_{Y}$:

\[ \int_{C}f(x,y) \ d_{\mu\times\mu'}(x,y)=\alpha(C)\]

\noindent As usual two measures are equal if and only if their associated 
functions agree almost everywhere.

Thus, we will define ${\cal N}(X,Y)$\ to be the set of all measures on
$X\times Y$ for which there exists a measurable function $f$ such that the
previous formula holds.  It is immediate that the marginals associated to
such a measure are absolutely continuous with respect to $\mu$\ and $\mu'$,
respectively.  While $f$ itself is only unique almost everywhere, the
measure with which $f$ is associated is easily viewed - in a canonical way -
both as a member of $Hom(X,Y)$ and as a member of $Hom(I,X \x Y)$.  Thus
every element of the set $Hom(I,X\ox Y)$ is associated with a measure that
has a functional kernel which is in turn one of the members of the set
${\cal N}(X,Y)$.  

To see that we have a 2-sided ideal, 
suppose that $\alpha\in{\cal N}(X,Y)$. Hence 
we have a function $f\colon X\times Y\rarr {\sf R}$ satisfying the above 
equation. Suppose $\beta\in Hom(Y,Z)$. let 
$G_2\colon \Sigma_Y\times Z\rarr {\sf R}$ be the associated
stochastic kernel. Then we define a function $h\colon X\times Z\rarr {\sf R}$
by the formula:

\[h(x,z)=\int_Y f(x,y)G_2(-,z)\]

As usual, we are viewing $f(x,y)$ as a measurable function of $y$ for the fixed
$x$, and $G_2(-,z)$ as a measure on $Y$ for the fixed $z$. The construction for
right composition is essentially identical. One can readily verify that the
functions so constructed are indeed functional kernels for
the composite measures. 

Finally, we observe that in the case when both $\alpha$ and $\beta$ are
nuclear, then there exist functions $f(x,y)$ and $g(y,z)$ which act as
functional kernels. The functional kernel of the composite is given by:

\[ \int_Y f(x,y)g(y,z)d\mu(y)\]

We conclude:

\begin{thm}
The above construction
determines a nuclear ideal for \prel{}. 
\end{thm}

The verification of the requirements for a nuclear ideal are routine.  The
calculations involve computing transposes and can be done just the same way
as proving associativity of composition in \cpd{}.  We call this nuclear
ideal \textbf{MRel}.  One can generalize the setting to \emph{analytic
spaces} \cite{Dudley89} which are continuous (or measurable) images of
$\mathbf{N}^{\infty}$ in Polish spaces.

\section{Trace Ideals}

\newcommand{\caL}{\mbox{${\cal L}$}}

In \cite{Joyal96}, Joyal, Street and Verity develop an abstract theory of
{\it trace operators} in a monoidal category. A {\it trace} 
is a function of the form:

\[ tr_A\colon Hom(A,A)\rarr Hom(I,I)\]

\noindent satisfying 
appropriate equations. (In fact, the authors introduce a more
general parametrized trace which we discuss below.) The authors
demonstrate that in a symmetric (in fact, braided) 
compact closed category, one obtains a trace
via the formula (using the notation of section 2 and using $c$ to represent 
the symmetry):

\[(h\colon A\rarr A)\mapsto (\nu;h\ox id;c;\psi\colon I\rarr I)\]

For example, in the compact closed category of finite-dimensional Hilbert
spaces, one obtains the usual notion of trace of an endomorphism. This notion
of trace also underlies such ideas as feedback in a computation and braid
closure \cite{Joyal96,Hasegawa97}. 

When one passes from the category of finite-dimensional Hilbert spaces 
to the category of arbitrary Hilbert spaces, one finds endomorphisms 
which do not have a trace, for example the identity on an 
infinite-dimensional space. However, each endomorphism monoid contains an ideal
of endomorphisms which do have a trace. This ideal is called 
the {\it trace class} and these trace maps are closely related to 
Hilbert-Schmidt morphisms. After reviewing this relationship, we describe
a general theory of {\it trace ideals} for symmetric monoidal categories. 
We then show that if a tensored $*$-category has a nuclear ideal satisfying
certain additional structure, then one can recover a trace ideal, as in the 
compact closed case.

\subsection{Hilbert Spaces}

Appropriate references for this material are \cite{Reed72,Simon79}.

\begin{defn} {\em 
An operator $B\in {\caL (H)}$, the space of bounded linear operators on
${\cal H}$, is called {\em positive} if $\langle Bx,x\rangle\geq 0$, for
all $x\in {\cal H}$. In this case, we write $B\geq 0$\ and $B\geq A$ if
$A-B\geq 0$.  }\end{defn}

Note for example that $AA^*$\ and $A^* A$ are always positive.

\begin{thm}(\cite{Reed72} page 196)
Suppose $A\geq 0$. Then there exists a unique $B\geq 0$ such that $B^{2}=A$. 
\end{thm}

\begin{defn} {\em The unique operator $B$ of the previous theorem is 
denoted $\sqrt{A}$.
Let $A\in {\caL (H)}$. Define $|A|=\sqrt{A^* A}.$}
\end{defn}

\begin{thm} Let ${\cal H}$ be separable
and  $\{e_i\}$ an orthonormal basis. If $A$ is
a positive operator, we define $tr(A)=\sum\langle Ae_n,e_n\rangle$. 
This is independent of
orthonormal basis. It has the following properties:

\begin{itemize}
\item $tr(A+B)=tr(A)+tr(B)$
\item $tr(\lambda A)=\lambda tr(A)$, for all $\lambda\geq 0$
\item If $0\leq A\leq B$, then $tr(A)\leq tr(B)$.
\end{itemize}
\end{thm}

\begin{defn}{\em
An operator $A$ is called {\it trace class} if $tr(|A|)<\infty$. The family 
of all trace class operators is denoted by ${\cal I}({\cal H})$ or just 
${\cal I}$.}
\end{defn}

\begin{thm}
${\cal I}$ has the following properties:

\begin{itemize}
\item ${\cal I}$ is a vector space.
\item It is a 2-sided ideal in the monoid $Hom({\cal H},{\cal H})$.
\item If $A\in {\cal I}$, then $A^*\in {\cal I}$
\end{itemize}
\end{thm}

These last two conditions say that we have a {\it $*$-ideal}.
We now extend the notion of trace to arbitrary endomorphisms in the trace 
ideal.

\begin{thm}(\cite{Reed72}, p.211)
If  $A\in {\cal I}$ and  $\{e_i\}$ is an orthonormal basis, then 
$\sum_{n=1}^{\infty}\langle Ae_n,e_n\rangle$ 
converges absolutely and is independent
of the basis. (We call this map the {\bf trace of $A$}, $tr(A)$.)
\end{thm}

Using the notion of trace class, it is possible to give an equivalent 
formulation of the notion of Hilbert-Schmidt  map:

\begin{prop}(\cite{Reed72}, p.211) 
A mapping $f\colon {\cal H}\rightarrow {\cal K}$ is
 Hilbert-Schmidt if and only if $f^*f\in
{\cal I}({\cal H})$. 
\end{prop}

The converse of this observation is also true:

\begin{prop}(\cite{Reed72}, p.211) 
If $h$ is a bounded linear operator on ${\cal H}$, then $h\in{\cal I}$ if 
and only if there exist Hilbert-Schmidt operators $f$ and $g$ on ${\cal H}$
such that
$h=fg$.
\end{prop}

\begin{rem}
Let ${\cal H}$ be a Hilbert space, and suppose 
we consider $H$  
as a Banach space. Then ${\cal H}$ is an object in the category
${\bf Ban}_{\infty}$, where we consider ${\bf Ban}_{\infty}$
with its usual $L_1$ tensor product. Thus we can apply Grothendieck's
original definition of nuclear morphism, and we see that we recover precisely 
the trace class maps.

\end{rem}

\subsection{Trace Ideals}

The previous discussion suggests the following abstract definition. We suppose
for the remainder that $\cal C$ is a symmetric monoidal category.

\begin{defn}{\em A {\em trace ideal} in $\cal C$ is a choice of subsets

\[ {\cal I}(U)\subseteq Hom(U,U) \mbox{  for each object $U$ in $\cal C$}\]

\noindent and a function

\[ tr_U\colon
{\cal I}(U)\rightarrow Hom(I,I) \mbox{  for each $U$ in $\cal C$}\]

\noindent such that

\begin{enumerate}

\item ${\cal I}(U)$ is a 2-sided ideal in the monoid $Hom(U,U)$.
\item ({\it Dinaturality} or {\it Sliding})
Suppose that $f\colon U\rarr V$ and $g\colon V\rarr U$ are such that 
$gf\in{\cal I}(U)$. Then $fg\in{\cal I}(V)$, and $tr_U(gf)=tr_V(fg)$.
\item ({\it Vanishing}) If $f\in {\cal I}(U)$, then $f\ox id_{I}\in
{\cal I}(U\ox I)$ and $tr_{U\ox I}(f\ox id_{I})=tr_{U}(f)$. Furthermore, 
we require that ${\cal I}(I)=Hom(I,I)$. If $f:I\rarr I$, then $tr_I(f)=f$.
\item ({\it Tensor Axiom}) If $f\in{\cal I}(U)$ and $g\in{\cal I}(V)$, then
$f\ox g\in {\cal I}(U\ox V)$ and $tr_{U\ox V}(f\ox g)=tr_U(f)tr_V(g)$.

\item Furthermore, if the category has a tensored $*$-structure, then we 
require that trace maps are closed under tensored $*$-structure, and
the trace operators respect this structure, i.e.

\begin{itemize}
\item If $f\in {\cal I}(U)$, then so is $f^*$, and $tr_{U}(f^*)=tr_{U}(f)^*$.
\item If $f\in {\cal I}(U)$, then $\overline{f}\in {\cal I}(\overline{U})$ and
$tr_{\overline{U}}(\overline{f})=tr_{U}(f)^*$.
\end{itemize}
\end{enumerate}
}\end{defn}

An alternative
approach to partial traces is presented in \cite{Blute98}, which 
considers traces on a {\it linearly distributive
category}. The trace operator works on a certain subcategory, the {\it core},
which has the same sort of ``type degeneracy'' as a compact closed category.

We would like to extend the relationship between compact closed
categories and traced monoidal categories to a relationship between
nuclear ideals and trace ideals. 
Keeping in mind the correspondence between Hilbert-Schmidt 
maps and the trace class, we define:

\begin{defn}{\em Suppose that $\cal C$ is a tensored $*$-category
equipped with a nuclear ideal.
Suppose also that $A$ is an object in $\cal C$. We define the {\em trace class}
of $A$ to be:}

\[ {\cal I}(A)=\{h\colon 
A\rightarrow A|\mbox{ There exists an object $B$, and morphisms
$f\colon A\rightarrow B,g\colon B\rarr A$}\] 
\[\mbox{with $f,g$ nuclear and $h=gf$}\}\]

{\em More generally, given two objects $A,B\in{\cal C}$, one can define: }

\[ {\cal I}(A,B)=\{h\colon 
A\rightarrow B|\mbox{ There exists an object $C$, and morphisms
$f\colon A\rightarrow C,g\colon C\rarr B$}\] 
\[\mbox{with $f,g$ nuclear and $h=gf$}\}\]
\end{defn}

\begin{lem}
${\cal I}(A)$ is a 2-sided ideal in the monoid $Hom(A,A)$.
${\cal I}(A,B)$ is a 2-sided ideal in $\cal C$.
\end{lem}

While one can define the notion of trace class for arbitrary 
morphisms in $\cal C$ as above, note that the actual trace function only acts
on ${\cal I}(A)={\cal I}(A,A)$. In other words, the trace function acts
only on the diagonal of the functor ${\cal I}(-,-)$. This is analogous
to the notion of {\it dinatural transformation}, which is the appropriate
notion of naturality for multivariate functors. These are families of 
morphisms between the two given functors, {\it instantiated along the 
diagonals}, satisfying an appropriate commutative hexagon 
\cite{EK,DS,BFSS,Blute93}. Hence the alternate name ``dinaturality''
for the sliding axiom.

If $h\in {\cal I}(A)$, we would like to define a morphism 
$tr_A(h)\colon I\rarr I$ 
(or just $tr(h)$ if there is no confusion) by the 
formula (where $\hat{g},\hat{f}$ denote the evident transposes):

\[ tr(h)=\hat{g}\hat{f}\colon I\rightarrow \overline{A}\ox B\rarr I\]

\noindent However, there is no guarantee that if $h$ is also equal to
$f'g'$ that we will obtain the same trace. Therefore we make the following 
definition:

\begin{defn}{\em
A nuclear ideal is {\em traced} if it satisfies the following uniqueness
property:

\begin{itemize}
\item If $f\colon A\rarr B, g\colon 
B\rarr A, f'\colon A\rarr C, g'\colon C\rarr A$ 
are nuclear and $gf=g'f'$, then 
$\hat{g}\hat{f}=\hat{g'}\hat{f'}\colon I\rarr I$.
\end{itemize}}
\end{defn}

\begin{thm}
The above construction assigns a trace ideal to each traced nuclear ideal.
\end{thm}

The proof of this theorem is simply a matter of checking 
the necessary diagrams. 
For example, lemma 5.8 gives the sliding axiom. One can
also check that:

\begin{thm}
The canonical nuclear ideal in ${\bf Hilb}$ is traced.
\end{thm}

\subsection{Traces in $\drel$}

We now examine the trace construction in our category of distributions.

\begin{thm} 
The canonical nuclear ideal in $\drel$ is traced.
\end{thm}

\begin{prf}{}
Suppose that $f\colon X\rarr Y, g\colon Y\rarr X, f'\colon X\rarr Z, g'\colon
Z\rarr X$ 
are nuclear and $gf=g'f'$. Since
$f$ is nuclear, we have a morphism $\hat{f}\colon I\rarr X\ox Y$, which has 
associated to it $\hat{f}_L\colon
\cD{I}\rarr \cD{X\times Y}$. As already remarked,
$\cD{I}$ is isomorphic to the base field, hence the map $\hat{f}_L$ simply
picks out an element of $\cD{X\times Y}$, which we denote by $\beta_f$. 
Similarly for $f',g,g'$. 

To verify the uniqueness property, recall that if $\phi\in\cD{X}$, then

\[ f_L(\phi)=\int_X\beta_f(x,y)\phi(x)\]

Since $gf=g'f'$, we have that for $\phi_1,\phi_2\in\cD{X}$:
 
\[\int_Yf_L(\phi_1)g_R(\phi_2)=\int_Zf'_L(\phi_1)g'_R(\phi_2) \]

After rearranging the order of integration one can conclude:

\[\int_X\int_X\int_Y\beta_f(x,y)\beta_g(y,x')\phi_1(x)\phi_2(x')=
\int_X\int_X\int_Z\beta_{f'}(x,z)\beta_{g'}(z,x')\phi_1(x)\phi_2(x')\]

\noindent The 
left-hand side corresponds to the distribution on $X\times X$ with kernel
$\int_Y\beta_f(x,y)\beta_g(y,x')$, and the right-hand side has kernel 
$\int_Z\beta_{f'}(x,z)\beta_{g'}(z,x')$. We know that two integrable functions
induce the same distribution if and only if they are equal almost everywhere,
but since these are smooth functions on 
$X\times X$, we conclude:

\[\int_Y\beta_f(x,y)\beta_g(y,x')=\int_Z\beta_{f'}(x,z)\beta_{g'}(z,x')\]

Thus we have:

\[\int_X\int_Y\beta_f(x,y)\beta_g(y,x)=
\int_X\int_Z\beta_{f'}(x,z)\beta_{g'}(z,x)\]

\noindent And we conclude $tr(gf)=tr(g'f')$.

\end{prf}

Actually, there is a more succinct description of the trace operator in \drel. 
Since $h=gf$ is nuclear, it has a kernel, $\alpha(x,x')$. Recall from theorem
\ref{drelnuc} that the formula for $\alpha$ is given by:

\[\alpha(x,x')=f_R(\beta_g(y,x'))=\int_Y \beta_f(x,y)\beta_g(y,x')\]

\noindent Hence we may conclude that:

\[tr_A(h)=\int_X\alpha(x,x)\]

We leave the details of the following to the reader. The result is quite 
similar to the case of \drel.

\begin{thm}
The canonical nuclear ideal in \prel \ is traced.
\end{thm}

\subsection{The parametric trace operator}

In \cite{Joyal96}, the authors actually have a parametrized trace operator. 
This means that there is a function of the form:

\[tr_U\colon Hom(A\ox U,B\ox U)\rarr Hom(A,B)\]

\noindent which reduces to the usual trace when $A=B=I$. There is an evident 
generalization to the ideal setting:

\begin{defn} {\em We suppose again 
that $\cal C$ is a symmetric monoidal category.
A {\em (parametric) trace ideal} 
in $\cal C$ is a choice of a family of 
subsets, for each object $U$ of $\cal C$ ,
of the form:

\[ {\cal I}^{U}_{A,B}\subseteq Hom(A\ox U,B\ox U) 
\mbox{  for all $A,B$ in $\cal C$}\]

\noindent and functions

\[ tr^{U}_{A,B}\colon {\cal I}^{U}_{A,B}\rightarrow Hom(A,B) \]

\noindent such that the families are ideals in the sense that:

\begin{itemize}
\item If $f\in{\cal I}^{U}_{A,B}$ and $h\colon U\rarr U$ is arbitrary, then 
$(id\ox h)\circ f$ and $f\circ (id\ox h)$ are in ${\cal I}^{U}_{A,B}$.
\item If $f\in{\cal I}^{U}_{A,B}$ and $g\colon B\rarr C, h\colon D\rarr A$
are arbitrary, then 
$(g\ox id_U)\circ f\circ (h\ox id_U)\in{\cal I}^{U}_{D,C}$. 
\end{itemize}

These are subject to the ideal-theoretic versions of the Joyal-Street-Verity 
axioms. In particular, (dropping sub- and superscripts if there is no
chance of confusion)

\begin{itemize}
\item (Vanishing) \begin{enumerate}
\item ${\cal I}^{I}_{A,B}=Hom(A\ox I,B\ox I)$, and the trace is calculated 
in the evident way.
\item Suppose $g\colon A\ox U\ox V\rarr B\ox U\ox V$. Then 
$g\in{\cal I}^{U\ox V}_{A,B}$ if and only if $g\in{\cal I}^{v}_{A\ox U,B\ox U}$
and $tr^{V}_{A\ox U,B\ox U}(g)\in{\cal I}^{U}_{A,B}$. Furthermore,

\[ tr^{U\ox V}_{A,B}(g)=tr^{U}_{A,B}(tr^{V}_{A\ox U,B\ox U}(g))\]
\end{enumerate} 
\item (Superposing) Suppose $f\in{\cal I}^{U}_{A,B}$ and $g\colon C\rarr D$ is
arbitrary. Then $g\ox f\in{\cal I}^{U}
_{C\ox A,D\ox B}$, and $tr(g\ox f)=g\ox tr(f)$.
\item (Yanking) 
Suppose $f\colon A\rarr U$ and $g\colon U\rarr B$. If 
$c_{U,B}\circ(f\ox g)\in {\cal I}^U_{A,B}$, then 
\[tr^U_{A,B}(c_{U,B}\circ(f\ox g))=gf\colon A\rarr B\]
\item (Sliding) Suppose $f\colon A\ox U\rarr B\ox V$ and $u\colon V\rarr U$. 
Then $(id\ox u)\circ f\in {\cal I}^{U}_{A,B}$ if and only if
$f\circ(id\ox u)\in {\cal I}^{V}_{A,B}$, and the two traces are equal.
\item (Tightening) Suppose 
$f\in{\cal I}^{U}_{A,B}$ and $g\colon B\rarr C, h\colon D\rarr A$
are arbitrary. Then 
\[tr((g\ox id_U)\circ f\circ (h\ox id_U))=g\circ tr(f)\circ h\]
\item Furthermore, if ${\cal C}$ is a tensored $*$-category, then 
the trace must preserve this structure in an evident sense.

\end{itemize}

}\end{defn}

Some discussion of our version of the Yanking axiom is in order. The 
Joyal-Street-Verity version of this axiom is essentially
the requirement that the trace of a symmetry morphism is the identity. However,
in our framework, one cannot make this requirement since the symmetry map will 
generally not be in the trace class. In the forthcoming thesis of Haghverdi
\cite{Haghverdi98}, it is observed that the following requirement is 
equivalent to the Joyal-Street-Verity version:

\bigskip

\noindent {\it Generalized Yanking Rule:}

Suppose $f\colon A\rarr U$ and $g\colon U\rarr B$. Then, 

\[tr^U_{A,B}(c_{U,B}\circ(f\ox g))=gf\colon A\rarr B\]

\subsection{$U$-nuclear ideals}

As before, we would like to construct trace ideals from nuclear ideals.
An analogous construction can be carried out using the notion of a 
{\it $U$-nuclear ideal}. 
We now outline this idea, but leave most of the details 
to the reader. The generalization amounts to introducing the notion of a {\it
$U$-nuclear morphism}. We will say that a morphism 
\mbox{$f\colon A\ox U\rarr B$} 
is $U$-nuclear, if it has a transpose $\hat{f}\colon A\rarr \overline{U}\ox B$.
More specifically, for each object $U$, we introduce a family of morphisms
${\cal N}_U(A\ox U,B)\subseteq Hom(A\ox U,B)$. These families should be
closed under all of the operations and furthermore
an ideal in the sense that if 

\[f\in{\cal N}_U(A\ox U,B)\] 

\noindent and $h\colon V\rarr U$ is arbitrary, then

\[ ((id\ox h);f)\in{\cal N}_V(A\ox V,B)\] 
 
\noindent Similarly for the variables $A$ and $B$.

Also there should be a natural bijection of the form:

\[ {\cal N}_V(A\ox V,B)\cong {\cal N}_{\overline{V}}(B\ox\overline{V},A)\]

\noindent satisfying appropriate equations. For example, the {\it compactness}
requirement becomes:

\begin{itemize}
\item (Compactness) Suppose $f\colon A\rarr C\ox B$ and $g\colon B\ox 
D\rarr E$. Then we have:

\end{itemize}

$$\begin{diagram}
A\ox D&\rTo^{f\ox id_D}&C\ox B\ox D\\
\dTo^{id_A \ox \hat{g} \ }&&\dTo_{ \ id_C \ox g}\\
A\ox\overline{B}\ox E&\rTo_{\hat{f}\ox id_E}&C\ox E\\
\end{diagram}$$

\smallskip

If a tensored $*$-category is equipped with such structure, we will refer to 
it as a {\it parametrized nuclear ideal.}

Given such a construction, one defines the {\it $U$-trace class} 
${\cal I}_U(A\ox U,B\ox U)\subseteq Hom(A\ox U,B\ox U)$ by saying that:

\[ h\in {\cal I}_U(A\ox U,B\ox U) \]

\noindent if and only if there exist 

\[f\in {\cal N}_U(A\ox U,C), g\in{\cal N}_U(B\ox U,C) \mbox{ \,\,\,\,\,  
such that } h=g^{*}f\]

One then constructs the $U$-trace of $h$ via the formula:

\[ tr^U_{A,B}(h)\colon A\rarr C\ox\overline{U}\rarr B\]

\noindent where the components are the evident transposes of $f$ and $h$. 
Again, one must add conditions to ensure that the trace satisfies appropriate 
equations. In particular, we note that with the above axioms, we can only
obtain the following weaker version of the yanking axiom:

\begin{lem} Suppose that $\cal C$ is
a tensored $*$-category equipped with a parametrized nuclear ideal. 
If $f\colon X\rarr U$ and $g\colon U\rarr Y$ are nuclear, then 
$c\circ (f\ox g)\colon X\ox U\rarr Y\ox U$ is in the $U$-trace class, and
\[ tr^U_{X,Y}(c\circ (f\ox g))=gf\]
\end{lem}

This is a consequence of the compactness requirement of section 5.

\subsection{Traces in \pmon}

We now discuss the traced structure of \pmon. First it is evident that unlike
in \hilb, we have that ${\cal I}(A)={\cal N}(A,A)$ for all objects $A$. If
$f\colon A\rarr A$ is a trace map, then we have the following formula:

\[tr(f)=\left\{ \begin{array}{ll}
id & \mbox{ if $|Dom(f)|=1$, and $f$ is the identity when 
restricted to its domain.}\\
\emptyset & \mbox{ otherwise}
\end{array}
\right. \]

The parametrized trace also has a very simple description. We will say that
a morphism $f\colon X\ox U\rarr Y$ is $U$-nuclear if it satisfies:

\[ \forall x\in X\mbox{ if $(x,u)\in Dom(f)$ and $(x,u')\in Dom(f)$, then 
$u=u'$}\]

Given this definition, there is an evident bijection
${\cal N}(X\ox U,Y)\cong{\cal N}(Y\ox U,X)$.

\noindent The class ${\cal I}(X\ox U,Y\ox U)$ is described by having the above 
requirement for both the domain and codomain. Then we can say that if
$f\in Tr(X\ox U,Y\ox U),~ (x,u)\in Dom(f)$ and $f(x,u)=(y,u')$, then:

\[tr(f)(x)\left\{\begin{array}{ll}
\mbox{undefined} & \mbox{ if $u\neq u'$}\\
y&\mbox{ if $u=u'$} 
\end{array}
\right. \]

\section{Conclusions}
Our investigations began with an attempt to define probabilistic
relations in analogy with ordinary relations.
Unexpectedly, ideas from functional analysis~\cite{Grothendieck55}
were essential.  The key idea, expressed in our abstract
definition of nuclear ideals, is that certain morphisms can be thought of
as behaving like ``matrices''.  

Our work naturally follows on from the development of Higgs and
Rowe~\cite{Higgs89}, the fundamental difference being that we have no
closed structure.  Crudely speaking, Higgs and Rowe generalize Banach
space theory while we generalize Hilbert space theory.

A key application of our work is that we can now work with structures
that are not categories but which are nuclear ideals inside some
tensored $*$-category.  For example, the nuclear ideal \mrel{}, described in 
Section 7,  
is of interest but is not a
category.  (As an example of its possible applications, we note that \mrel{}
has partially additive structure ~\cite{Manes86,Haghverdi98}.) 
However, \mrel{} is indeed a
nuclear ideal in \prel{}.

An important open question is the computational significance of trace ideals.
It is already well-established that a trace structure can be used to model
feedback in denotational semantics \cite{Joyal96,Hasegawa97}.  But what can be
said when one only has these operations on an ideal? The geometry of 
interaction program, due to Girard \cite{Girard88}, can be used to obtain 
a compact closed category from a traced monoidal category
\cite{Abramsky94a,Abramsky96,Joyal96}.  It seems possible that a similar
construction applied to a category with a traced ideal will give a nuclear 
ideal.  

Another area of application of the theory of compact closed categories is
{\em topological quantum field theory} \cite{Atiyah89,Atiyah90}, which
evolved, in part, from Segal's work on {\em conformal field theory} 
\cite{Segal88}.  In topological quantum field theory, one considers a
compact closed category of {\em cobordisms} in which composition is defined
by gluing along boundaries.  Then a TQFT is given by a compact closed functor 
to the compact closed category of finite-dimensional Hilbert spaces.  
In Segal's formulation of conformal field theory, one works with arbitrary
Hilbert spaces and a similar ``category'' of Riemann surfaces with 
boundary.  This structure is essentially a compact closed category, 
except that it fails to be a category in that it lacks identity 
morphisms.  Thus it seems reasonable to 
suspect that it is a nuclear ideal in some larger ambient tensored 
$*$-category.  One of our goals in future work will be to find such a
category.  A conformal field theory would then be a {\em nuclear functor}
to the tensored $*$-category {\bf Hilb}.

A related issue is the extension of our work to higher-dimensional
categories.  The theory of {\em $n$-Hilbert spaces} \cite{Baez96}, a
higher-dimensional analogue of Hilbert space, 
has become quite important in TQFT \cite{Baez95}.  Baez has developed
the theory of $2$-Hilbert spaces with this in mind, and extended some of the 
work of Doplicher and Roberts to this setting \cite{Doplicher89}.  

Finally, the category \drel{} suggests several further topics of investigation.
One possible extension of \drel{} is to the theory of {\em noncommutative 
distributions} \cite{Albeverio93}.  Roughly speaking, these are distributions
which take values in a Lie group.  They are useful in the representation
theory of gauge groups.  Finally, we hope to take advantage of the fact that
distributions form a {\em $\cal D$-module}, 
that is to say they provide representations of the Weyl algebra
\cite{Coutinho95}.  It would be interesting to attempt to extend the work
of \cite{Blute96,Blute97}, where full completeness theorems are obtained by
considering representations of the additive group of integers and a 
noncocommutative Hopf algebra.

\bibliography{nuclear}

\begin{thebibliography}{10}

\bibitem{Abramsky94}
S.~Abramsky, S.~Gay, and R.~Nagarajan.
\newblock Interaction categories and foundations of typed concurrent
  programming.
\newblock In M.~Broy, editor, {\em Deductive Program Design: Proceedings of the
  1994 Marktoberdorf International Summer School}, NATO ASI Series F.
  Springer-Verlag, 1994.
\newblock Also available as theory/papers/Abramsky/marktoberdorf.ps.gz via
  anonymous ftp to theory.doc.ic.ac.uk.

\bibitem{Abramsky94a}
S. Abramsky, R. Jagadeesan.
\newblock New foundations for the geometry of interaction.
\newblock {\em Information and Computation}, 111(1):53--119, May 1994.

\bibitem{Abramsky96}
S.~Abramsky.
\newblock Retracing some paths in process algebra.
\newblock In Montanari and Sassone, editors, {\em Proceedings of CONCUR 96},
  number 1119 in Lecture Notes In Computer Science, pages 1--17.
  Springer-Verlag, 1996.

\bibitem{Ageron96}
P.~Ageron.
\newblock Effective taxonomies and crossed taxonomies.
\newblock {\em Cahiers de Top. et Geom. Diff.}, 37:82--90, 1996.

\bibitem{Albeverio93}
S.~Albeverio, R.~Hoegh-Krohn, J.~Marion, D. Testard, and B. Torr\'esani, 
\newblock {\it Noncommutative Distributions}.
\newblock Dekker Pure and Applied Mathematics, 1993.

\bibitem{Algwaiz92}
M.~A. Al-Gwaiz.
\newblock {\em Theory of Distributions}.
\newblock Dekker Pure and Applied Mathematics, 1992.

\bibitem{Ash72}
R.~B. Ash.
\newblock {\em Real Analysis and Probability}.
\newblock Academic Press, 1972.

\bibitem{Atiyah89}
M.~Atiyah.
\newblock Topological quantum field theories.
\newblock {\em Publ. Math. Inst. Hautes Etudes Sci. Paris}, 68:175--186, 1989.

\bibitem{Atiyah90}
M.~Atiyah.
\newblock {\em The Geometry and Physics of Knots}.
\newblock Cambridge University Press, 1990.

\bibitem{Baez95}
J.~Baez, J.~Dolan.
\newblock Higher-dimensional algebra and topological quantum field theory.
\newblock {\em Journal of Mathematical Physics,} 36:6073-6105, 1995.

\bibitem{Baez96}
J.~Baez.
\newblock Higher-dimensional algebra II:$2$-Hilbert spaces.
\newblock preprint, 1996.

\bibitem{BFSS}
E. Bainbridge, P. Freyd, A. Scedrov, P. Scott.
\newblock Functorial Polymorphism
\newblock {\it Theoretical Computer Science,} 70:35-64, 1990.

\bibitem{Barr80}
M.~Barr.
\newblock {\em $*$-autonomous categories}.
\newblock Number 752 in Lecture Notes in Mathematics. Springer-Verlag, 1980.

\bibitem{Billingsley95}
P.~Billingsley.
\newblock {\em Probability and Measure}.
\newblock Wiley-Interscience, 1995.

\bibitem{Blute93}
R.~Blute.
\newblock Linear logic, coherence and dinaturality.
\newblock {\em Theoretical Computer Science}, 115:3-41, 1993.

\bibitem{Blute96}
R.~Blute, P.~Scott.
\newblock Linear L\"auchli semantics.
\newblock {\em Annals of Pure and Applied Logic}, 77:101-142, 1996.

\bibitem{Blute97}
R.~Blute, P.~Scott.
\newblock The shuffle Hopf algebra and noncommutative full completeness.
\newblock to appear in {\em Journal of Symbolic Logic,} 1998.

\bibitem{Blute98}
R. Blute, J.R.B. Cockett, and R.A.G. Seely. 
\newblock Feedback for linearly distributive categories: traces and fixed 
points.
\newblock In preparation, 1998.

\bibitem{Carboni87}
A.~Carboni and R.~F.~C. Walters.
\newblock Cartesian bicategories i.
\newblock {\em Journal of Pure and Applied Algebra}, 49:11--32, 1987.

\bibitem{Coutinho95}
S.~Coutinho.
\newblock {\em A primer of algebraic $\cal D$-modules}.
\newblock London Mathematical Society Student Texts,
Cambridge University Press, 1995.

\bibitem{Danos}
V. Danos, 
\newblock {\it Logique lin\'eaire: Une repr\'esentation alg\'ebrique du 
calcul,}
\newblock preprint


\bibitem{Dieudonne88}
J.~Dieudonne.
\newblock {\em Treatise on Analysis - {VII}}.
\newblock Number~10 in Pure and Applied Mathematics. Academic Press, 1988.


\bibitem{Doplicher89}
S.~Doplicher and J.~Roberts.
\newblock A new duality theory for compact groups.
\newblock {\em Inventiones Mathematicae}, 98:157--218, 1989.

\bibitem{DS}
E. Dubuc, R. Street.
\newblock Dinatural Transformations
\newblock {\em Springer Lecture Notes in Mathematics Volume 137}, 
Springer-Verlag, 1970

\bibitem{Dudley89}
R.~M. Dudley.
\newblock {\em Real Analysis and Probability}.
\newblock Wadsworth and Brookes/Cole, 1989.

\bibitem{EK}
S. Eilenberg, G.M. Kelly. 
\newblock A generalization of the functorial calculus
\newblock {\em Journal of Algebra}, 3:366-375, 1966.

\bibitem{Freyd90}
P.~J. Freyd, A.~Scedrov.
\newblock {\em Categories, Allegories}.
\newblock North-Holland, 1990.

\bibitem{Freyd89}
P. Freyd, D. Yetter.
\newblock Braided compact closed categories with applications to low 
dimensional topology.
\newblock {\em Advances in Mathematics}, 77:156-182, 1989.


\bibitem{Ghez85}
P. ~Ghez, R.~Lima, and J.~Roberts.
\newblock $w^*$-categories.
\newblock {\em Pacific Journal of Mathematics}, 120:79--109, 1985.

\bibitem{Girard87}
J.-Y. Girard.
\newblock Linear logic.
\newblock {\em Theoretical Computer Science}, 50:1--102, 1987.

\bibitem{Girard88}
J.Y. Girard.
\newblock Geometry of interaction I: interpretation of system $F$. 
\newblock {\em Proceedings of the ASL Meeting, Padova, 1988.}

\bibitem{Giry80}
M.~Giry.
\newblock A categorical approach to probability theory.
\newblock In B.~Banaschewski, editor, {\em Proceedings of a Conference on
  Categorical Aspects of Topology and Analysis}, number 915 in Lecture Notes In
  Mathematics, pages 68--85. Springer-Verlag, 1980.

\bibitem{Grothendieck55}
A.~Grothendieck.
\newblock {\em Products Tensoriels Topologiques et Espaces Nucleaires}.
\newblock AMS Memoirs. American Mathematical Society, 1955.

\bibitem{Halmos50}
P.~Halmos.
\newblock {\em Measure Theory}.
\newblock Graduate Texts in Mathematics 18. Springer-Verlag, 1974.
\newblock Originally published in 1950.

\bibitem{Haghverdi98}
E. Haghverdi. 
\newblock Thesis, In preparation, 1998. 

\bibitem{Hasegawa97}
M. Hasegawa.
\newblock Recursion from cyclic sharing: traced monoidal categories and 
models of cyclic lambda calculi.
\newblock {\it Springer Lecture Notes in Computer Science 1210}, 
p.196-213, 1997.

\bibitem{Higgs89}
D.~A. Higgs, K.~Rowe.
\newblock Nuclearity in the category of complete semilattices.
\newblock {\em Journal of Pure and Applied Algebra}, 57:67--78, 1989.

\bibitem{Hoffman-Jorgenson94}
J. Hoffman-J\"orgenson.
\newblock {\em Probability With a View Towards Applications - 2 volumes.}
\newblock Chapman and Hall, 1994.


\bibitem{Hormander90}
L. H\"ormander.
\newblock {\em The Analysis of Linear Partial Differential Operators I.}
\newblock Grundleheren der mathematischen Wissenscaften 256, Springer-Verlag 
1990.

\bibitem{Joyal96}
A.~Joyal, R.~Street and D.~Verity.
\newblock Traced monoidal categories.
\newblock {\em Mathematical Proceedings of the Cambridge Philosophical 
Society,} 119:425--446, 1996.

\bibitem{Joyal84}
A.~Joyal, M.~Tierney.
\newblock {\em An Extension of the Galois Theory of Grothendieck}.
\newblock Memoirs of the AMS. American Mathematical Society, 1984.



\bibitem{Kadison83}
R. ~Kadison, J.~Ringrose.
\newblock {\em Fundamentals of the Theory of Operator Algebras}.
\newblock Academic Press, 1983.

\bibitem{MacLane71}
S. Mac Lane.
\newblock {\em Categories for the Working Mathematician}, volume~5 of {\em
  Graduate texts in Mathematics}.
\newblock Springer-Verlag, New York, 1971.

\bibitem{Malliavin95}
P.~Malliavin.
\newblock {\em Integration and Probability}.
\newblock Graduate Texts in Mathematics 157. Springer-Verlag, 1995.
\newblock French edition appeared in 1993.

\bibitem{Manes86}
E.~Manes, M.~Arbib.
\newblock {\em Algebraic Approaches to Program Semantics}.
\newblock Springer-Verlag, 1986.

\bibitem{Pachl78}
J. Pachl.
\newblock Disintegration and Compact Measures.
\newblock {\em Math. Scand.}, 43:157-168, 1978

\bibitem{Porter95}
T. Porter. 
\newblock Interpretations of Yetter's notion of $G$-coloring: simplicial fibre
bundles and nonabelian cohomology.
\newblock preprint, 1995.

\bibitem{Raney60}
G.~Raney.
\newblock Tight Galois connections and complete distributivity.
\newblock {\em Transactions of the American Mathematical Society,} 97:418--426,
(1960)

\bibitem{Reed72}
M. Reed, B. Simon.
\newblock {\em Functional Analysis, Methods of Mathematical Physics, Volume I}.
\newblock Academic Press 1972

\bibitem{Rowe88}
K.~A. Rowe.
\newblock Nuclearity.
\newblock {\em Canad. Math. Bull.}, 31(2):227--235, 1988.

\bibitem{Rudin66}
W.~Rudin.
\newblock {\em Real and Complex Analysis}.
\newblock McGraw-Hill, 1966.

\bibitem{Schwartz57}
L.~Schwartz.
\newblock {\em Th\'eorie des Distributions}.
\newblock Hermann, 1957.


\bibitem{Segal88}
G.~Segal.
\newblock The definition of conformal field theory.
\newblock In K.~Bleuler and M.~Werner, editors, {\em Differential Geometric
  Methods in Theoretical Physics}, pages 165--171. Kluwer Academic Publishers,
  1988.

\bibitem{Selinger97}
P.~Selinger.
\newblock First order axioms for concurrency.
\newblock In {\em Proceedings of CONCUR 97}, number 1243 in Lecture Notes In
  Computer Science, 1997.

\bibitem{Simon79}
B. Simon, 
\newblock {\em Trace Ideals and Their Applications}
\newblock Cambridge University Press, 1979

\bibitem{Treves67}
F.~Treves.
\newblock {\em Topological Vector Spaces, Distributions and Kernels}.
\newblock Pure and Applied Mathematics 25. Academic Press, 1967.


\bibitem{Wendt93}
M.~Wendt.
\newblock {\em On Measurably Indexed Families of Hilbert Spaces}.
\newblock PhD thesis, Dalhousie University, 1993.

\bibitem{Wendt94}
M.~Wendt.
\newblock The category of disintegrations.
\newblock {\em Cahiers de Topologie et Geometrie Differentielle Categoriques},
  35:291--308, 1994.

\bibitem{Yetter}
D. Yetter.
\newblock Topological quantum field theories associated to finite groups and
crossed $G$-sets.
\newblock {\em Jourbal of Knot Theory and its Ramifications} 2:113-123, 1993.

\end{thebibliography}

\end{document}